\documentclass[11pt]{article}
\usepackage{srcltx,graphicx}
\usepackage{amsmath, amssymb, amsthm}
\usepackage{color, xcolor, colortbl}
\usepackage{array, multirow}
\usepackage{subfig, overpic}
\usepackage{hyperref}
\usepackage{epstopdf}
\usepackage{algorithm}
\usepackage{algorithmic}
\RequirePackage[capitalize,nameinlink]{cleveref}[0.19]
\usepackage{ulem}

\graphicspath{{images/}}

\newtheorem{property}{Property}

\newcommand\bbR{\mathbb{R}}

\newcommand\bbN{\mathbb{N}}

\newcommand\cH{\mathcal{H}}
\newcommand\cI{\mathcal{I}}
\newcommand\cJ{\mathcal{J}}
\newcommand\cM{\mathcal{M}}
\newcommand\cP{\mathcal{P}}
\newcommand\cT{\mathcal{T}}

\newcommand\NL{\mathsf{NL}}
\newcommand\IL{\mathsf{IL}}
\newcommand\LCR{\mathsf{LCR}}
\newcommand\LCK{\mathsf{LCK}}
\newcommand\LCI{\mathsf{LCI}}
\newcommand\CNNR{\mathsf{CR}}
\newcommand\CNNK{\mathsf{CK}}
\newcommand\CNNI{\mathsf{CI}}
\newcommand\Reshape{\mathsf{Reshape}}
\newcommand\Flatten{\mathsf{Flatten}}
\newcommand\linear{\mathsf{linear}}

\newcommand\Ntrainsample{$N^{\mathrm{train}}_{\mathrm{samples}}$}
\newcommand\Ntestsample{$N^{\mathrm{test}}_{\mathrm{samples}}$}
\newcommand\Nparams{$N_{\mathrm{params}}$}

\newcommand\ad{\mathrm{ad}}
\newcommand\ie{{\it i.e.}~}

\newcommand\sps[1]{^{(#1)}}

\newcommand{\ud}{\,\text{d}}

\newcommand\dd{\,\mathrm{d}}     
\numberwithin{equation}{section}

\usepackage[top=1in, bottom=1in, left=1in, right=1in]{geometry}

\newcommand\revised[2]{#2} 
\newcommand\add[1]{\revised{}{#1}}

\title{A multiscale neural network based on hierarchical matrices}
\author{
Yuwei Fan\thanks{Department of Mathematics, Stanford University,
    Stanford, CA 94305, email: {\tt ywfan@stanford.edu}},~~
Lin Lin\thanks{Department of Mathematics, University of California, Berkeley, and Computational
    Research Division, Lawrence Berkeley National Laboratory, Berkeley, CA 94720,
    email: {\tt linlin@math.berkeley.edu}},~~
Lexing Ying\thanks{Department of Mathematics and Institute for
    Computational and Mathematical Engineering, Stanford University,
    Stanford, CA 94305, email: {\tt lexing@stanford.edu}},~~
Leonardo Zepeda-N\'u\~nez\thanks{Computational
    Research Division, Lawrence Berkeley National Laboratory, Berkeley, CA 94720 , email: {\tt lzepeda@lbl.gov }}
}

\begin{document}
\date{}
\maketitle
\begin{abstract}
In this work we introduce a new multiscale artificial neural network based on the structure of
$\mathcal{H}$-matrices. This network generalizes the latter to the nonlinear case by introducing a
local deep neural network at each spatial scale. Numerical results indicate that the network is able
to efficiently approximate discrete nonlinear maps obtained from discretized nonlinear partial
differential equations, such as those arising from nonlinear Schr\"odinger equations and the
Kohn-Sham density functional theory.
\end{abstract}

\vspace*{4mm}

\section{Introduction}\label{sec:intro}
In the past decades, there has been a great combined effort in developing efficient algorithms to
solve linear problems issued from discretization of integral equations (IEs), and partial
differential equations (PDEs). In particular, multiscale methods such as multi-grid methods
\cite{Brandt1977}, the fast multipole method \cite{greengard1987fast}, wavelets \cite{Mallat2009},
and hierarchical matrices \cite{borm2003introduction,hackbusch1999sparse}, have been strikingly
successful in reducing the complexity for solving such systems. In several cases, 
for operators of pseudo-differential type, these algorithms can achieve linear or
quasi-linear complexity. In a nutshell, these methods aim
to use the inherent multiscale structure of the underlying physical problem to build efficient
representations at each scale, thus compressing the information contained in the system. The gains
in complexity stem mainly from processing information at each scale, and merging it in a
hierarchical fashion.

Even though these techniques have been extensively applied to linear problems with outstanding
success, their application to nonlinear problems has been, to the best of our knowledge, very
limited. This is due to the high complexity of the solution maps.
In particular, building a global approximation of such maps would
normally require an extremely large amount of parameters, which in return, is often translated 
to algorithms with a prohibitive computational cost. The development of algorithms and heuristics 
to reduce the cost is an area of active research
\cite{Barrault2004,EFENDIEV2007,efendiev_galvis_li_presho_2014,GreplMaday2007,MadayMula2016}.
However, in general, each method is application-dependent, and requires a deep understanding of the underlying physics.

On the other hand, the surge of interest in machine learning methods, in particular, deep neural
networks, has dramatically improved speech recognition \cite{Hinton2012}, visual object recognition
\cite{Krizhevsky2012}, object detection, etc. This has fueled breakthroughs in many domains such as
drug discovery \cite{MaSheridan2015}, genomics \cite{Leung2014}, and automatic translation
\cite{SutskeverNIPS2014}, among many others \cite{leCunn2015,SCHMIDHUBER2015}. Deep neural networks have empirically shown that it is possible to obtain
efficient representations of very high-dimensional functions.  Even though the universality theorem
holds for neural networks \cite{CohenSharir2018,Hornik91,Khrulkov2018,Mhaskar2018}, {\it i.e.}, they can
approximate arbitrarily well any function with mild regularity conditions, how to efficiently build such approximations remains an open question. In particular, the degree of approximation 
depends dramatically on the architecture of the neural network, \ie how the different layers are 
connected. In addition, there is a fine balance between the number of parameters, the architecture, 
and the degree of approximation \cite{He2015ConvolutionalNN,He2016DeepRL,Mhaskar2018}.

This paper aims to combine the tools developed in deep neural networks with ideas from multiscale
methods.  In particular, we extend hierarchical matrices ($\cH$-matrices) to nonlinear problems 
within the framework of neural networks. Let
\begin{equation}\label{eq:non_linear_map}
  u = \cM(v),\quad u,v \in \Omega \subset \mathbb{R}^n,
\end{equation}
\add{be a nonlinear generalization of pseudo-differential operators,}
issued from an underlying physical problem, described by
either an integral equation or a partial differential equation, where $v$ can be considered as a parameter in the equation, $u$ is either the
solution of the equation or a function of it, and $n$ is the number of variables.
  
We build a neural network with a novel multiscale structure inspired by hierarchical matrices.
We interpret the application of an $\cH$-matrix to a vector using a neural network structure as follows.
We first reduce the dimensionality of the vector, or {\it restrict} it, by multiplying it by a short and wide structured 
matrix. Then we {\it process} the encoded vector by multiplying it by a structured square matrix. Then we
return the vector to its original size, or {\it interpolate} it, by multiplying it by a structured
tall and skinny matrix. 
Such operations are performed separately at different spatial scales. 
The final approximation to the matrix-vector multiplication is obtained by adding the
contributions from all spatial scales, including the near-field contribution, which is represented by a near-diagonal matrix.
Since every step is linear, the overall operation is also a linear mapping. 
This interpretation allows us to directly generalize the $\cH$-matrix to
nonlinear problems by replacing the structured square matrix in the processing stage by a structured nonlinear network with several layers. The resulting artificial neural network, which we call \emph{multiscale neural network}, only
requires a relatively modest amount of parameters even for large problems.  

We demonstrate the performance of the multiscale neural network by approximating the solution to the
nonlinear Schr\"odinger equation \cite{anglin2002bose,pitaevskii1961vortex}, as well as the
Kohn-Sham map \cite{HohenbergKohn1964, KohnSham1965}. These mappings are highly nonlinear, and are still well approximated by the proposed neural network, with a relative accuracy on the order of $10^{-4}\sim
10^{-3}$. Furthermore, we do not observe overfitting even in the presence of a relatively small
set of training samples. 

\subsection{Related work}
Although machine learning and deep learning literature is vast, the application of deep
learning to numerical analysis problems is relatively new, though that is 
rapidly changing.
Research using deep neural networks with multiscale architectures has primarily focused on
image \cite{badrinarayanan2015segnet,Bruna2012,Chen2018DeepLab,LITJENS201760,Ronneberger2015,Ulyanov2018} and video
\cite{Mathieu2016} processing.

Deep neural networks have been recently used to solve PDEs
\cite{Beck2017,berg2017unified,Chaudhari2017,E2017,han2018solving,Yingzhou2018,Raissi2018DL,Rudd2014,Konstantinos2018} and 
classical inverse problems \cite{Araya-Polo2018,PASCHALIS2004211}.  
For general applications of machine learning to nonlinear numerical analysis problems, the work of
Raissi and Karnidiakis used machine learning, in particular, Gaussian processes, to find parameters
in nonlinear equations \cite{Raissi2018}; 
Chan, and Elsheikh predicted the basis function on the coarse grid in multiscale finite volume
method by neural network;
Khoo, Lu and Ying used neural network in the context of uncertainty quantification
\cite{khoo2017solving}; 
Zhang et al used neural network in the context of generating high-quality interatomic potentials for
molecular dynamics \cite{ZhangHanWangCarE2018,ZhangWangE2018}. 
Wang et al. applied non-local multi-continuum neural network on time-dependent nonlinear problems
\cite{WangChung2018}.
Khrulkov et al. and Cohen et al. developed deep neural network architectures based on tensor-train
decomposition \cite{CohenSharir2018,Khrulkov2018}. 
\add{In addition, we note that deep neural networks with
related multi-scale structures \cite{Ronneberger2015,PeltSethian2018} have been proposed mainly for applications such as image processing, however, 
we are not aware of any applications of such architectures to solving nonlinear differential or integral equations.}

\subsection{Organization}

The reminder of the paper is organized as follows. \Cref{sec:NNhmatrix} reviews the
$\cH$-matrices and interprets the $\cH$-matrices using the framework of neural
networks. \Cref{sec:hnn} extends the neural network representation of $\cH$-matrices to the
nonlinear case. \Cref{sec:application} discusses the implementation details and
demonstrates the accuracy of the architecture in representing nonlinear maps,
followed by the conclusion and future directions in \cref{sec:conclusion}.

 \section{Neural network architecture for $\cH$-matrices}\label{sec:NNhmatrix}

In this section, we aim to represent the matrix-vector multiplication of $\cH$-matrices within the
framework of neural networks. For the sake of clarity, we succinctly review the structure of
$\cH$-matrices for the one dimensional case in \cref{sec:Hmatreview}. We interpret
$\cH$-matrices using the framework of neural networks in 
\cref{sec:NNrepresentation}, and then extend it to the multi-dimensional case in 
\cref{sec:nD}.

\subsection{$\cH$-matrices} \label{sec:Hmatreview}
Hierarchical matrices ($\cH$-matrices) were first introduced by Hackbusch et al. in a series of
papers \cite{hackbusch1999sparse, hackbusch2000sparse,hackbusch2001introduction} as an algebraic
framework for representing matrices with a hierarchically off-diagonal low-rank structure. This framework provides
efficient numerical methods for solving linear systems arising from integral equations (IE) and partial
differential equations (PDE) \cite{borm2003introduction}. In the sequel, we follow the notation in
\cite{lin2011fast} to provide a brief introduction to the framework of $\cH$-matrices in a simple
uniform and Cartesian setting. The interested readers are referred to
\cite{hackbusch1999sparse,borm2003introduction,lin2011fast} for further details.

Consider the integral equation
\begin{equation}\label{eq:integral}
    u(x) = \int_{\Omega}g(x,y)v(y) \dd y, \quad \Omega=[0,1),
\end{equation}
where $u$ and $v$ are periodic in $\Omega$ and $g(x,y)$ is smooth and numerically low-rank away from
the diagonal. A discretization with an uniform grid with $N=2^Lm$ discretization points yields the
linear system given by
\begin{equation}\label{eq:discrete}
    u = A v,
\end{equation}
where $ A\in\bbR^{N\times N}$, and
   $u, v \in\bbR^N$ are the discrete analogues of $u(x)$ and $v(x)$ respectively.

We introduce a hierarchical dyadic decomposition of the grid in $L+1$ levels as follows.
We start by the $0$-th level of the decomposition, which corresponds to the set of all grid points defined as
\begin{equation}\label{eq:cI}
    \cI\sps{0} = \{k/N: k=0,\dots,N-1\}.
\end{equation}
At each level $\ell$ ($0\leq \ell\leq L$), we decompose the grid in $2^{\ell}$ disjoint {\em
segments}.

Each segment is defined by $\cI\sps{\ell}_i = \cI\sps{0} \cap[(i-1)/2^{\ell}, i/2^{\ell})$ for
$i=1,\dots,2^{\ell}$. Throughout this manuscript, $\cI\sps{\ell}$(or $\cJ\sps{\ell}$) denote
a generic segment of a given level $\ell$, and the superscript $\ell$ will be omitted when the level
is clear from the context. Moreover, following the usual terminology in $\cH$-matrices, we say that a
segment $\cJ\sps{l}$ ($\ell\geq1$) is the \emph{parent} of a segment $\cI\sps{l-1}$ if $\cI\sps{l-1}
\subset \cJ\sps{l}$. Symmetrically, $\cI\sps{l-1}$ is called a \emph{child} of $\cJ\sps{l}$.
Clearly, each segment, except those on level $L$, have two child segments.

\begin{figure}[ht]
  \centering
  \includegraphics[width=\textwidth]{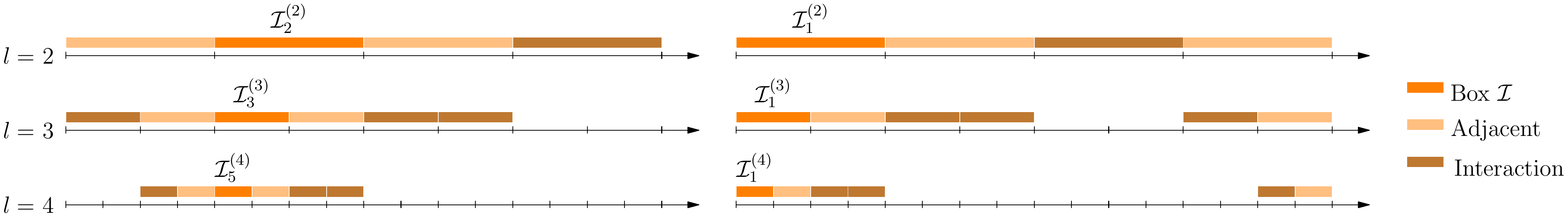}
  \caption{\label{fig:partitionmesh}Illustration of the computational domain at level $\ell=2,3,4$.
    The left and right figures represent an interior segment and a boundary segment and their adjacent and
    interaction list at different levels.}
\end{figure}
In addition, for a segment $\cI$ on level $\ell\geq 2$, we define the following lists:
\begin{itemize}
\item[$\NL(\cI)$] \emph{neighbor list} of  $\cI$. List of the segments on level $\ell$ that are
  adjacent to $\cI$ including $\cI$ itself;
\item[$\IL(\cI)$] \emph{interaction list} of $\cI$.  If $\ell = 2$, $\IL(\cI)$ contains all the
  segments on level $2$ minus $\NL(\cI)$. If $\ell > 2$, $\IL(\cI)$ contains all the segments on
  level $\ell$ that are children of segments in $\NL(\cP)$ with $\cP$ being $\cI$'s parent minus
  $\NL(\cI)$.
\end{itemize}

\cref{fig:partitionmesh} illustrates the dyadic partition of the computational domain and the
lists on levels $\ell=2,3,4$.  Clearly, $\cJ\in\NL(\cI)$ if and only if $\cI\in\NL(\cJ)$, and
$\cJ\in\IL(\cI)$ if and only if $\cI\in\IL(\cJ)$.

\begin{figure}[ht]
    \centering
    \includegraphics[width=\textwidth]{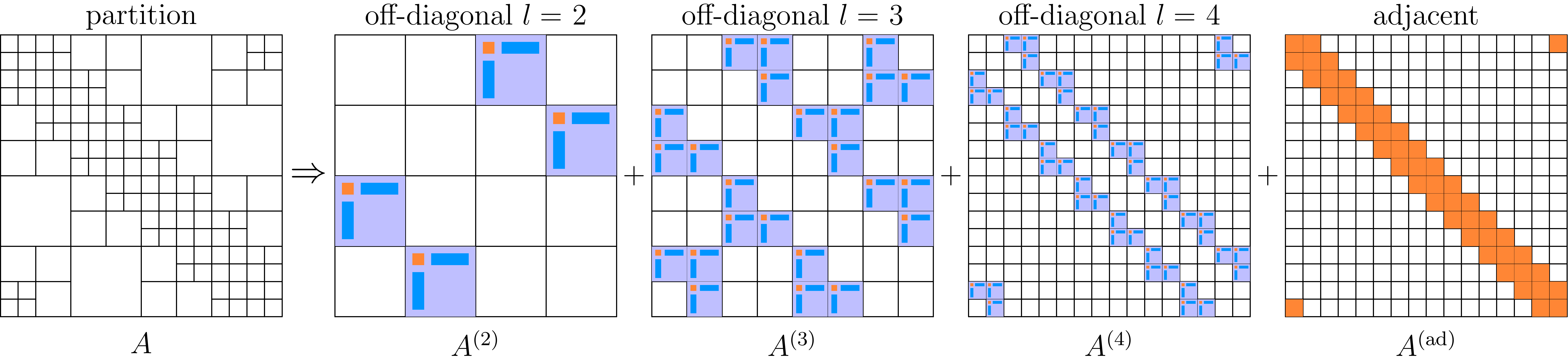}
    \caption{\label{fig:partitionmatrix}Partition of the matrix $A$ for $L=4$ and nonzero blocks of
    $A\sps{\ell}$, $\ell=2,3,4$ (colored blue) and $A^{(\ad)}$ (colored orange). Nonzero blocks of
    $A\sps{\ell}$ are approximated by low rank approximation and nonzero blocks of $A^{(\ad)}$ are
    retained.}
\end{figure}

For a vector $v\in\bbR^N$, $v_{\cI}$ denotes the elements of $v$ indexed by $\cI$ and for a matrix
$A\in\bbR^{N\times N}$, $A_{\cI,\cJ}$ represents the submatrix given by the elements of $A$ indexed
by $\cI\times \cJ$. The dyadic partition of the grid and the interaction lists induce a multilevel
decomposition of $A$ as follows
\begin{equation}\label{eq:decompose}
    A = \sum_{\ell=2}^L A\sps{\ell} + A^{(\ad)},
    \qquad
    \begin{aligned}
    A\sps{\ell}_{\cI,\cJ} &= \begin{cases}
        A_{\cI,\cJ}, & \cI \in \IL(\cJ);\\
        0, & \text{otherwise},
    \end{cases}
    \cI,\cJ \text{ at level } \ell, 2\le \ell \le L,\\
    A^{(\ad)}_{\cI,\cJ} &= \begin{cases}
        A_{\cI,\cJ}, & \cI \in \NL(\cJ); \\
        0, & \text{otherwise},
    \end{cases}
    \cI,\cJ \text{ at level } L. 
    \end{aligned}
\end{equation}

In a nutshell, $A\sps{\ell}$ considers the interaction at level $\ell$ between a
segment and its interaction list, and $A^{(\ad)}$ considers all the interactions between adjacent
segments at level $L$.

\cref{fig:partitionmatrix} illustrates the block partition of $A$ induced by the dyadic
partition, and the decomposition induced by the different interaction lists at each level that
follows \eqref{eq:decompose}.

\begin{figure}[ht]
    \centering
    \subfloat[$u\sps{\ell}=A\sps{\ell}v\approx U\sps{\ell}M\sps{\ell}(V\sps{\ell})^Tv$]{
    \label{fig:Avl}
    \begin{overpic}[height=0.15\textheight]{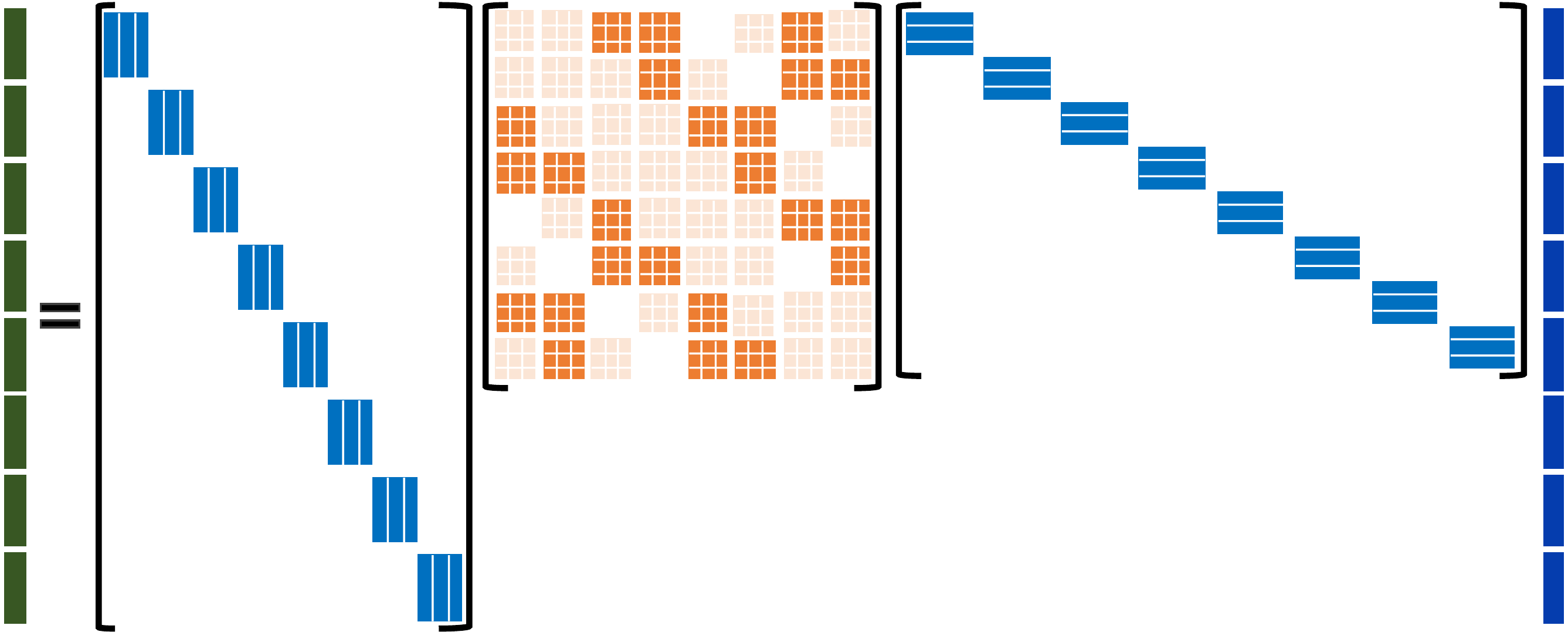}
        \put(-1,42){$u\sps{\ell}$}
        \put(15,42){$U\sps{\ell}$}
        \put(40,42){$M\sps{\ell}$}
        \put(70,42){$(V\sps{\ell})^T$}
        \put(98,42){$v$}
    \end{overpic}
    }\quad    \rule{0.5pt}{0.15\textheight}    \quad    \subfloat[$u^{(\ad)}=A^{(\ad)}v$]{
    \label{fig:Avad}
    \begin{overpic}[height=0.15\textheight]{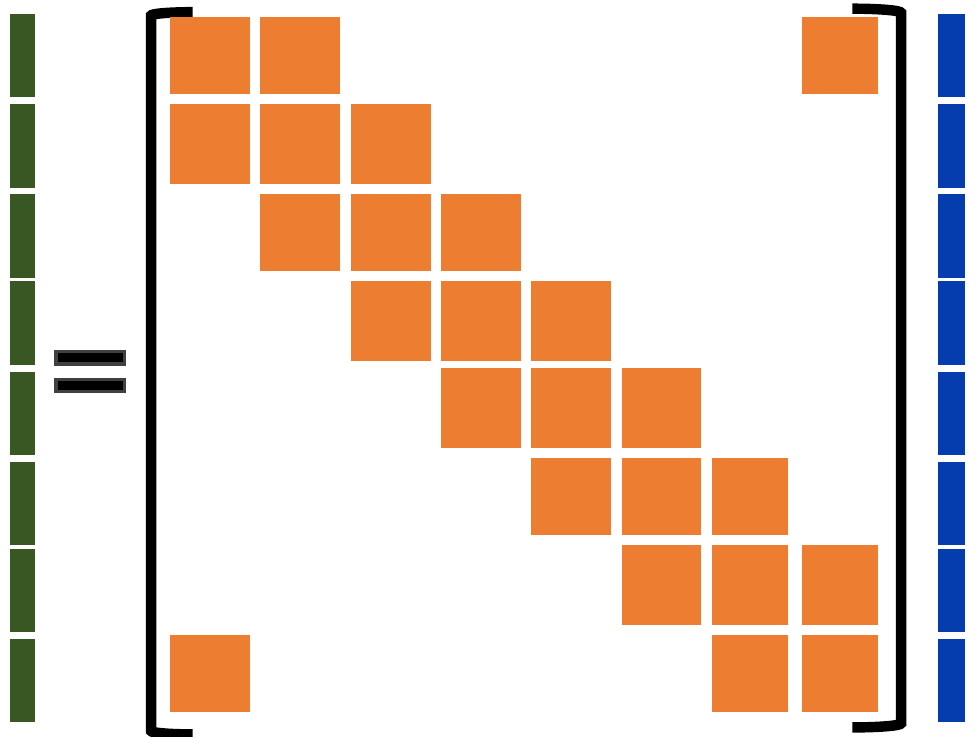}
        \put(-1,78){$u^{(\ad)}$}
        \put(50,78){$A^{(\ad)}$}
        \put(95,78){$v$}
    \end{overpic}
    }
    \caption{\label{fig:Av}Diagram of matrix-vector multiplication \eqref{eq:H-Av} of the
      low-rank part and adjacent part of $\cH$-matrices. The blocks of $M\sps{\ell}$ colored by
      pale orange are zero blocks, and if we treat these blocks as non-zero blocks, the matrices
      $M\sps{\ell}$ are block cyclic band matrices.}
\end{figure}
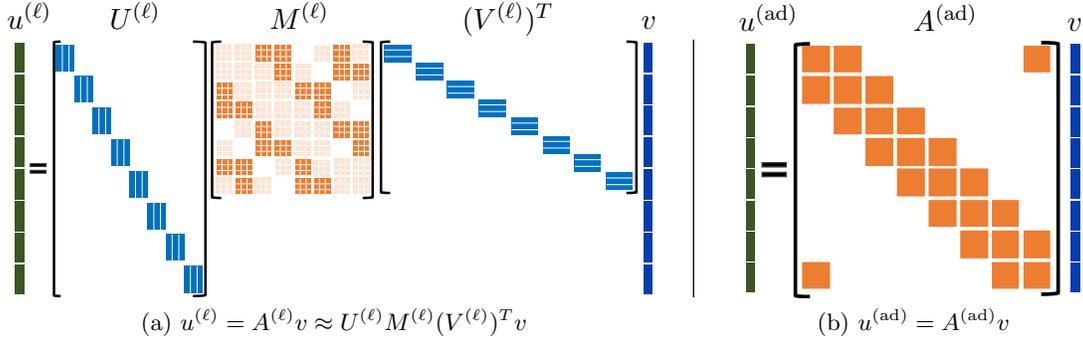

The key idea behind $\cH$-matrices is to approximate the nonzero blocks $A\sps{\ell}$ by a low rank
approximation (see \cite{kishore2017literature} for a thorough review). 
This idea is depicted in
\cref{fig:partitionmatrix}, in which each non-zero block is approximated by a tall and skinny
matrix, a small squared one and a short and wide one, respectively. In this work, we focus on the uniform
$\cH$-matrices \cite{fenn2002fmm}, and, for simplicity, we suppose that each block has a fixed rank
at most $r$, \ie
\begin{equation} \label{eq:low_rank_fact}
  A\sps{\ell}_{\cI, \cJ}\approx U\sps{\ell}_{\cI}M\sps{\ell}_{\cI,\cJ}(V\sps{\ell}_{\cJ})^{T},
    \quad U\sps{\ell}_{\cI}, V\sps{\ell}_{\cJ}\in\bbR^{N / 2^{\ell}\times r},\quad
    M\sps{\ell}_{\cI,\cJ}\in\bbR^{r\times r}.
\end{equation}
where $\cI,$ and $ \cJ$ are any interacting segments at level $\ell$.

The main observation is that it is possible to factorize each $A\sps{\ell}$ as $A\sps{\ell} \approx
U\sps{\ell}M\sps{\ell}(V\sps{\ell})^T$ depicted in \cref{fig:Av}. $U\sps{\ell}$ is a block
diagonal matrix with diagonal blocks $U\sps{\ell}_{\cI}$ for $\cI$ at level $\ell$, $V\sps{\ell}$ is
a block diagonal matrix with diagonal blocks $V\sps{\ell}_{\cJ}$ for $\cJ$ at level $\ell$, and
finally $M\sps{\ell}$ aggregates all the blocks $M\sps{\ell}_{\cI,\cJ}$ for all interacting segments
$\cI,\cJ$ at level $\ell$. This factorization induces a decomposition of $A$ given by
\begin{equation}\label{eq:Hmatrix}
  A = \sum_{\ell=2}^LA\sps{\ell}+A^{(\ad)}
  \approx \sum_{\ell=2}^LU\sps{\ell}M\sps{\ell}(V\sps{\ell})^T+A^{(\ad)}.
\end{equation}
Thus the matrix-vector multiplication \eqref{eq:discrete} can be expressed as
\begin{equation}\label{eq:H-Av}
  \begin{aligned}
    u & \approx  \sum_{\ell=2}^L u\sps{\ell} + u^{(\ad)} = \sum_{\ell=2}^L
    U\sps{\ell}M\sps{\ell}(V\sps{\ell})^T v + A^{(\ad)} v,
    \end{aligned}
\end{equation}
as illustrated in \cref{fig:Av}, which constitutes the basis for writing $\cH$-matrices as a neural network.

In addition, the matrices $\{U\sps{\ell},V\sps{\ell},M\sps{\ell}\}_{\ell=2}^{L}$ and $A^{(\ad)}$ have the
following properties.
\begin{property}\label{pro:A}
    The matrices
    \begin{enumerate}
        \item $U\sps{\ell}$ and $V\sps{\ell}$, $\ell=2,\cdots,L$ are block diagonal matrices with block size $N /
            2^{\ell}\times r$; \label{pro:UV}
        \item $A^{(\ad)}$ is a block cyclic tridiagonal matrix with block size $m\times m$; \label{pro:Aad}
        \item $M\sps{\ell}$, $\ell=2,\cdots,L$ are block cyclic band matrices with block size $r\times r$
            and band size $n_{b}\sps{\ell}$, which is $2$ for $\ell=2$ and $3$ for $\ell\geq3$, if we treat all
            the  pale orange colored blocks of $M\sps{\ell}$ in \cref{fig:Avl} as nonzero blocks.
            \label{pro:M}
    \end{enumerate}
\end{property}
We point out that the band sizes $n_{b}\sps{\ell}$ and $n^{(\ad)}_{b}$
depend on the definitions of $\NL$ and $\IL$. In this case, the list were defined using the {\em strong
admissible condition} in $\cH$-matrices. Other conditions can be certainly used, such as the {\em
weak admissibility condition}, leading to similar structures.

\subsection{Matrix-vector multiplication as a neural network} \label{sec:NNrepresentation}

An artificial neural network, in particular, a feed-forward network, can be thought of the 
composition of several simple functions, usually called \emph{propagation functions},
in which the intermediate one-dimensional variables are called \emph{neurons},  which in return, are
organized in vector, or tensor, variables called \emph{layers}. For example, an artificial
feed-forward neural network 
\begin{equation}
    u =  \mathcal{F}(v), \qquad u, v \in \bbR^n
\end{equation}
with $K+1$ layer can be written using the following recursive formula 
\begin{equation} \label{eq:neuralnet}
\begin{split}
    \xi\sps{0} &= v, \\
\xi\sps{k} &= \phi(W\sps{k}\xi\sps{k-1} + b\sps{k}),\\
u &= \xi\sps{K},
\end{split}
\end{equation}
where for each $k = 1,...,K$ 
we have that $ \xi\sps{k},b\sps{k} \in \bbR^{n_{k}},
W\sps{k} \in \bbR^{n_{k}\times n_{k-1}}$.  Following the terminology of machine
learning, $\phi$ is called the \emph{activation function} that is applied entry-wise, $W\sps{k}$
are the \emph{weights}, $b\sps{k}$ are the \emph{biases}, and $\xi\sps{k}$ is the $k$-th layer
containing $n_{k}$ number of \emph{neurons}. Typical choices for the activation function are linear
function, the rectified-linear unit (ReLU), or sigmoid function. In addition, \eqref{eq:neuralnet}
can easily be rewritten using tensors by replacing the matrix-vector multiplication by the more
general tensor contraction. We point out that representing layers with tensors or vectors is
equivalent up to reordering and reshaping. The main advantage of using the former is that layers,
and its propagating functions, can be represented in a more compact fashion. Therefore, in what 
follows we predominantly use a tensor representation. 

Within this context, training a network refers to finding the weights and biases, whose entries
are collectively called \emph{parameters},  in order to
approximate a given map. This is usually done by minimizing a loss function using a stochastic
optimization algorithm.

\subsubsection{Locally connected networks}
We interpret the structure of $\cH$-matrices \eqref{eq:Hmatrix} using the framework of neural
networks.  The different factors in \eqref{eq:H-Av} possess a distinctive
structure, which we aim to exploit by using locally connected (LC) network. LC networks are
propagating functions whose weights have a block-banded constraint. For the one-dimensional example, 
we also treat $\xi$ as a 2-tensor of dimensions $\alpha \times N_x$, where $\alpha$ is the {\em channel
  dimension} and $N_x$ is the {\em spatial dimension}, and $\zeta$ be a 2-tensor of dimensions
$\alpha' \times N_x'$. We say that $\xi$ is connected to $\zeta$ by a LC networks if
\begin{equation}\label{eq:lc}
    \zeta_{c',i} = \phi\left( \sum_{j=(i-1)s+1}^{(i-1)s+w}\sum_{c=1}^{\alpha}
    W_{c',c;i,j}\xi_{c,j} + b_{c',i}\right),
    \quad i=1,\dots,N_{x}', ~ c'=1,\dots,\alpha',
\end{equation}
where $w$ and $s\in\bbN$ are the {\em kernel window size} and {\em stride}, respectively. 
In addition, we say that $\zeta$ is a {\em locally connected} (LC) layer if it satisfies \eqref{eq:lc}.

\begin{figure}[ht]
    \begin{minipage}{\textwidth}
        \begin{tabular}{ >{\centering\arraybackslash}m{0.31\textwidth}|
            >{\centering\arraybackslash}m{0.30\textwidth}|
            >{\centering\arraybackslash}m{0.30\textwidth}}
                                                                                    \includegraphics[width=0.31\textwidth]{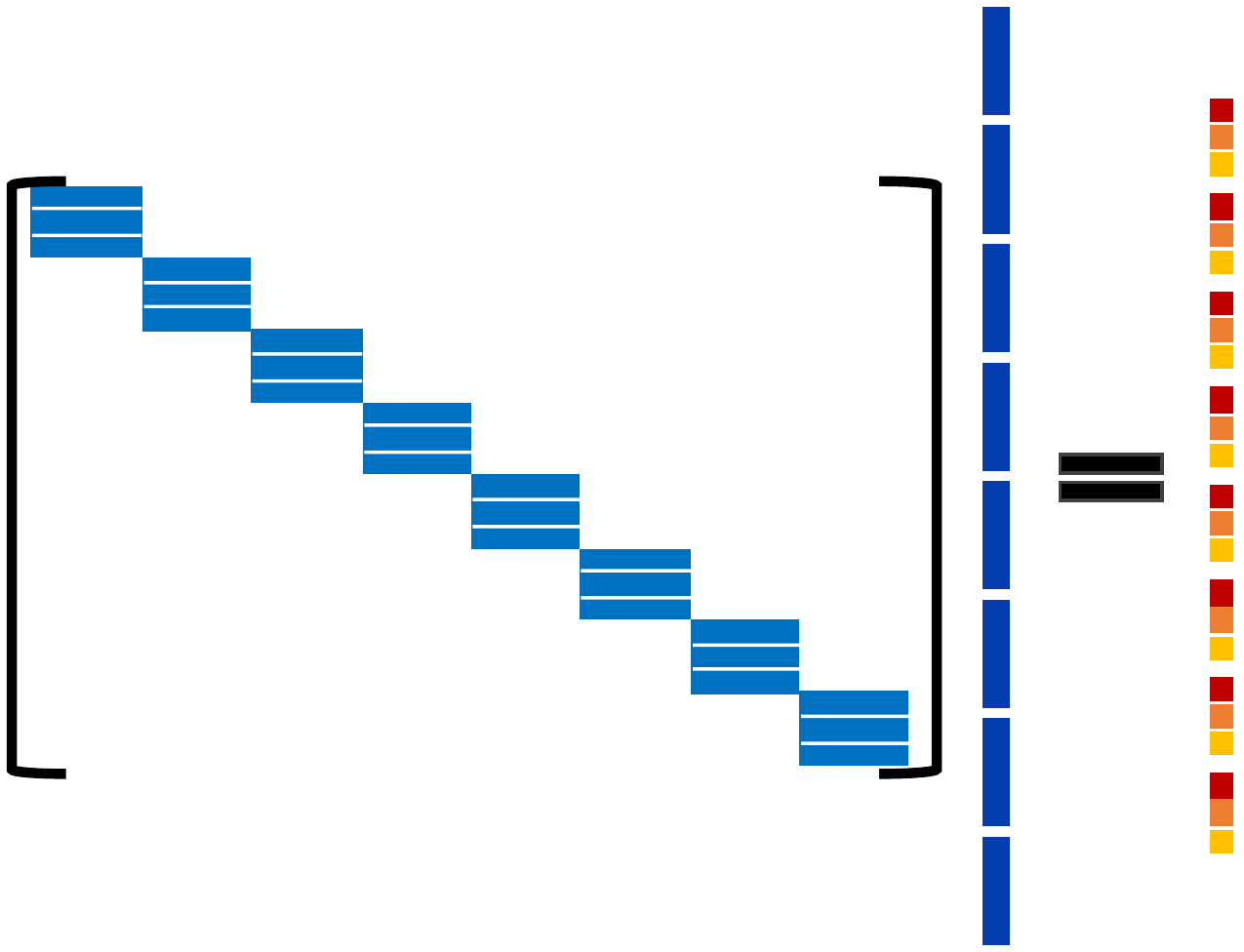} &
            \includegraphics[width=0.3\textwidth]{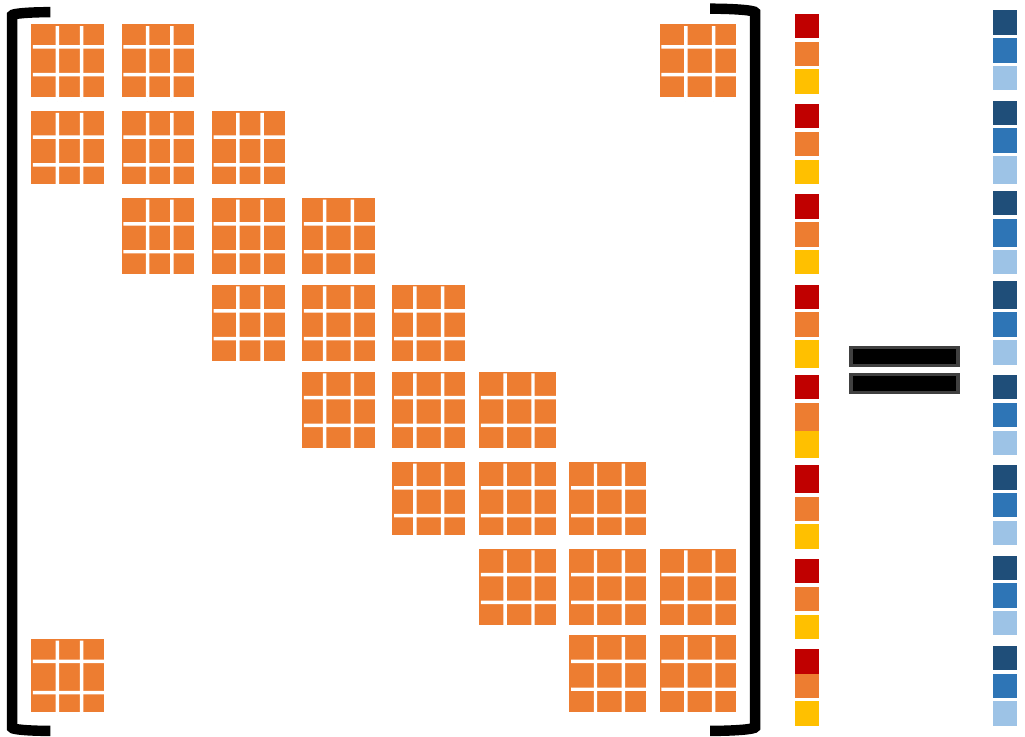} &
            \includegraphics[width=0.23\textwidth]{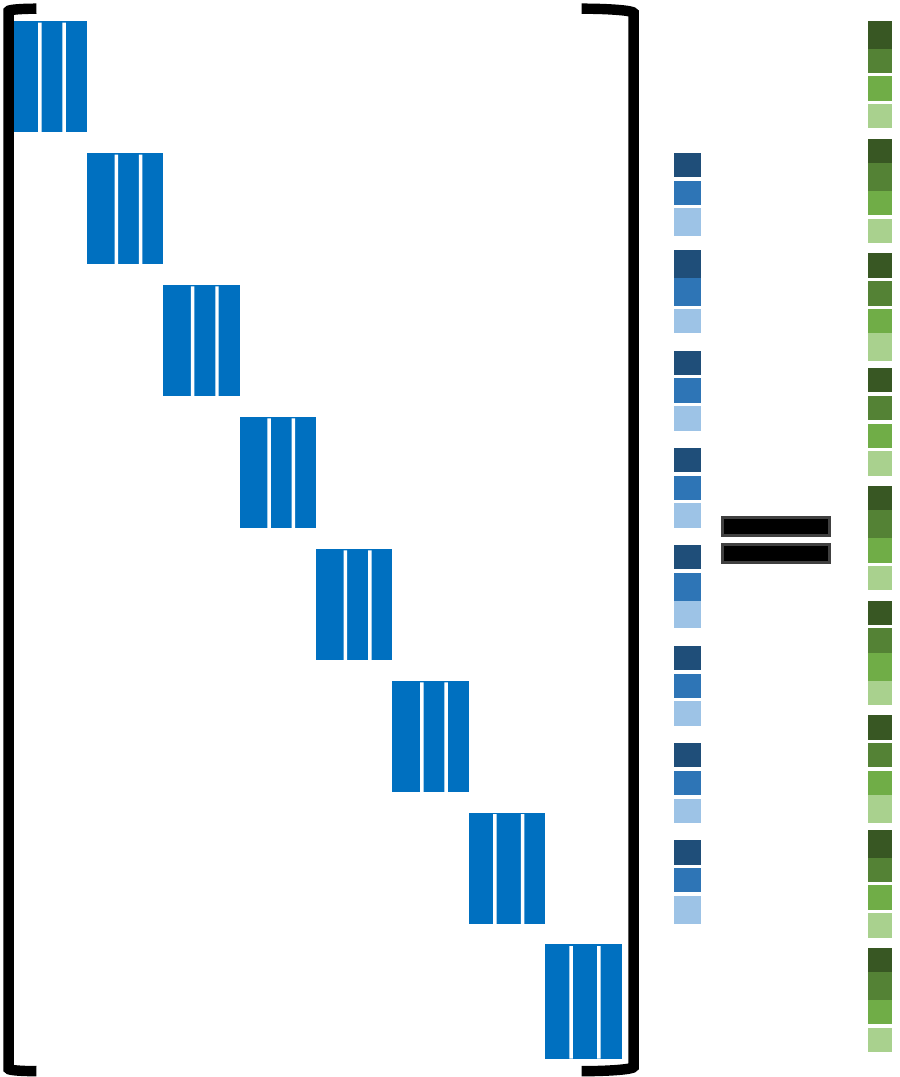} \\[2mm]
            \includegraphics[page=1, width=0.30\textwidth]{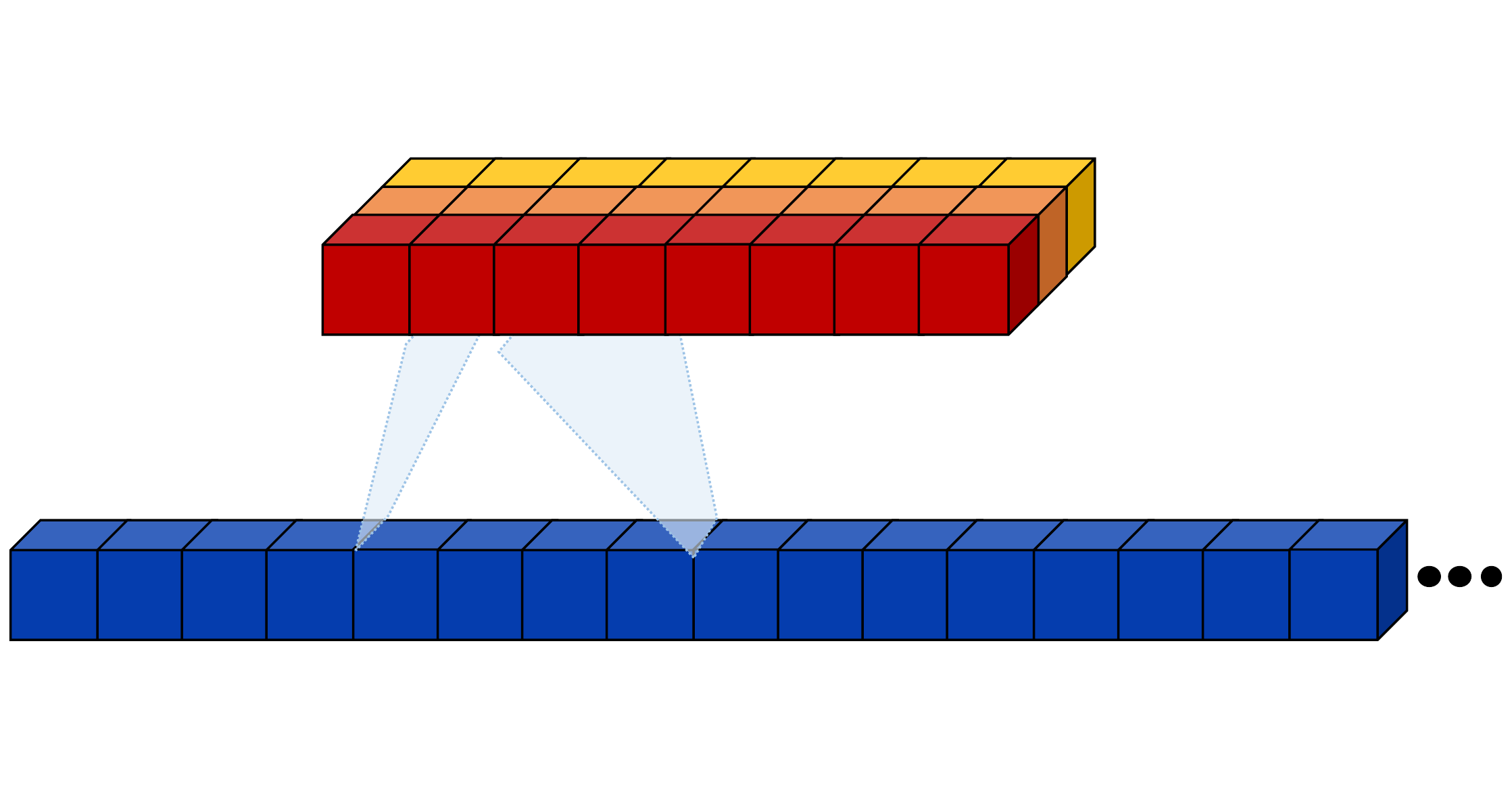} &
            \includegraphics[page=2, width=0.30\textwidth]{LCs.pdf} &
            \includegraphics[page=3, width=0.30\textwidth]{LCs.pdf} \\
            $s=w=\frac{N_x}{N_x'}$, $\alpha=1$  &
            $s=1$, $N_x'=N_x$ &
            $s=1$, $w=1$, $N_x'=N_x$ \\
            (a) $\LCR[\phi; N_x, N_x', \alpha']$ with 
            $N_x=32$, $N_x'=8$ and $\alpha'=3$ &
            (b) $\LCK[\phi; N_x, \alpha, \alpha', w]$ 
            with $N_x=8$, $\alpha=\alpha'=3$ and $w=3$ &
            (c) $\LCI[\phi; N_x, \alpha, \alpha']$ 
            with $N_x=8$, $\alpha=3$ and $\alpha'=4$ \\
        \end{tabular}
    \end{minipage}
    \caption{\label{fig:LC}Three instances of locally connected networks used the represent the matrix-vector multiplication in \eqref{eq:H-Av}. The upper portions of each column depicts the patterns of the matrices and the lower portions are their respective analogues using locally connect networks.}
\end{figure}

Each LC network requires $6$ parameters, $N_x$, $\alpha$, $N_x'$, $\alpha'$, $w$ and $s$ to be characterized. Next, we
define three types of LC network by specifying some of their parameters,
\begin{itemize}
    \item[$\LCR$] \emph{Restriction} network: we set $s=w=\frac{N_x}{N_x'}$ and $\alpha=1$ in LC.
        This network represents the multiplication of a block diagonal matrix with block sizes
        $\alpha'\times s$ and a vector with size $N_x\alpha$, as illustrated by \cref{fig:LC} (a). 
        We denote this network using $\LCR[\phi; N_x, N_x', \alpha']$.
        The application of $\LCR[\linear; 32, 8, 3]$ is depicted in \cref{fig:LC} (a). 

    \item[$\LCK$] 
        \emph{Kernel} network: we set $s=1$ and $N_x'=N_x$.  This network represents the multiplication of a
        cyclic block band matrix of block size $\alpha'\times \alpha$ and band size
        $\frac{w-1}{2}$ times a vector of size $N_x\alpha$, as illustrated by the upper portion of
        \cref{fig:LC} (b). To account for the periodicity we pad the input layer $\xi_{c,j}$ on
        the spatial dimension to the size $(N_x+w-1)\times \alpha$. We denote this network by
        $\LCK[\phi; N_x, \alpha, \alpha', w]$. This network has two steps: the periodic padding of
        $\xi_{c,j}$ on the spatial dimension, and the application of \eqref{eq:lc}. 
        The application of $\LCK[\linear; 8,3, 3, 3]$ is depicted in \cref{fig:LC} (b).

    \item[$\LCI$] \emph{Interpolation} network: we set $s=w=1$ and $N_x'=N_x$ in LC. This network
        represents the multiplication of a block diagonal matrix with block size $\alpha'\times
        \alpha$, times a vector of size $N_x\alpha$, as illustrated by the upper figure in 
        \cref{fig:LC} (c). We denote the network $\LCI[\phi; N_x, \alpha, \alpha']$, which has two steps:
        the application of \eqref{eq:lc}, and the reshaping of the output to a vector by column
        major indexing.
        The application of $\LCI[\linear; 8, 3, 4]$ is depicted in \cref{fig:LC} (c).
\end{itemize}

\subsubsection{Neural network representation}

Following \eqref{eq:H-Av}, in order to construct the neural network (NN) architecture for \eqref{eq:H-Av}, we need
to represent the following four operations:
\begin{subequations}\label{eq:hm}
    \begin{align}
        \xi\sps{\ell} & =(V\sps{\ell})^Tv, \label{eq:hmI}\\
        \zeta\sps{\ell} &=M\sps{\ell}\xi\sps{\ell},\label{eq:hmK}\\
        u\sps{\ell}&=U\sps{\ell}\zeta\sps{\ell},\label{eq:hmR}\\
        u^{(\ad)}&=A^{(\ad)}v. \label{eq:hmad}
    \end{align}
\end{subequations}
From \cref{pro:A}.\ref{pro:UV} and the definition of $\LCR$ and $\LCI$, the
operations \eqref{eq:hmI} and \eqref{eq:hmR} are equivalent to
\begin{equation}\label{eq:lcR}
    \xi\sps{\ell} = \LCR[\linear; N,2^{\ell},r](v),\quad
    u\sps{\ell} = \LCI[\linear; 2^{\ell}, r, \frac{N}{2^{\ell}}](\zeta\sps{\ell}),
\end{equation}
respectively.
Analogously, \cref{pro:A}.\ref{pro:M} indicates that \eqref{eq:hmK} is equivalent to
\begin{equation}\label{eq:lcM}
    \zeta\sps{\ell} = \LCK[\linear; 2^{\ell},r, r, 2n_b\sps{\ell}+1](\xi\sps{\ell}).
\end{equation}
We point out that $\xi\sps{\ell}$  is a vector in \eqref{eq:hmR} but a $2$-tensor in
\eqref{eq:lcR} and \eqref{eq:lcM}.
In principle, we need to flatten $\xi\sps{\ell}$ in \eqref{eq:lcR} to a vector and reshape it back to a
2-tensor before \eqref{eq:lcM}. These operations do not alter the algorithmic pipeline, so they are omitted.

Given that $v,u^{(\ad)}$ are vectors, but $\LCK$ is defined for 2-tensors, we
explicitly write the reshape and flatten operations. Denote as $\Reshape[n_1,n_2]$ the map that reshapes a
vector of size $n_1n_2$ into a 2-tensor of size $n_1\times n_2$ by column major
indexing, and $\Flatten$ is defined as the inverse of $\Reshape$. 
Using \cref{pro:A}.\ref{pro:Aad}, we can write \eqref{eq:hmad} as
\begin{equation}\label{eq:lcAd}
        \tilde{v} = \Reshape[m,2^L](v), ~
    \tilde{u}^{(\ad)}=\LCK\left[\linear; 2^L,m,m,2n^{(\ad)}_{b}+1\right](\tilde{v}), ~
    u^{(\ad)}=\Flatten(\tilde{u}^{(\ad)}).
\end{equation}

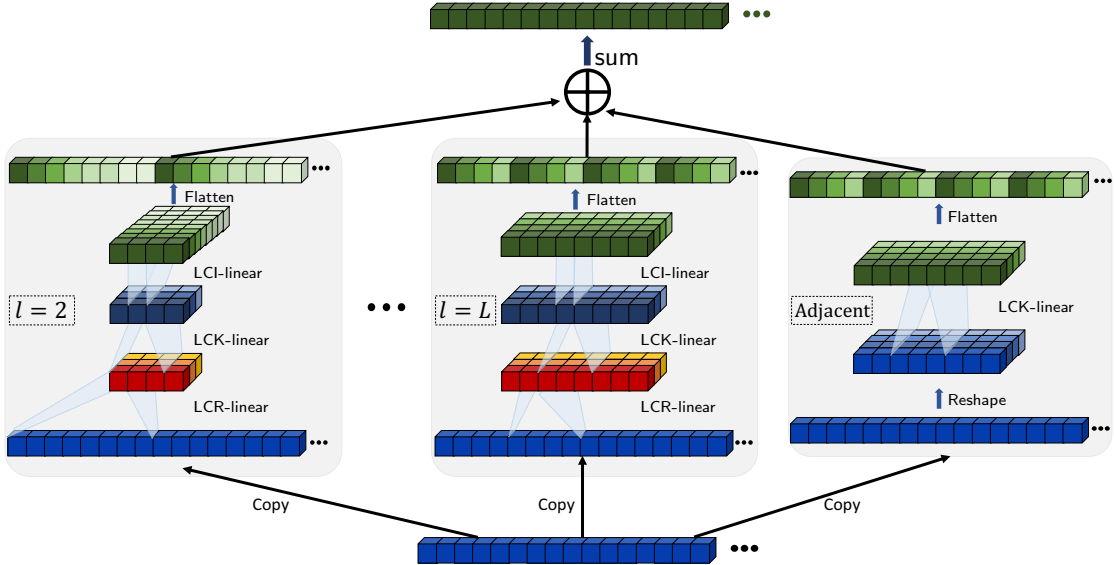
\begin{figure}[htb]
    \centering
    \begin{overpic}[width=0.9\textwidth]{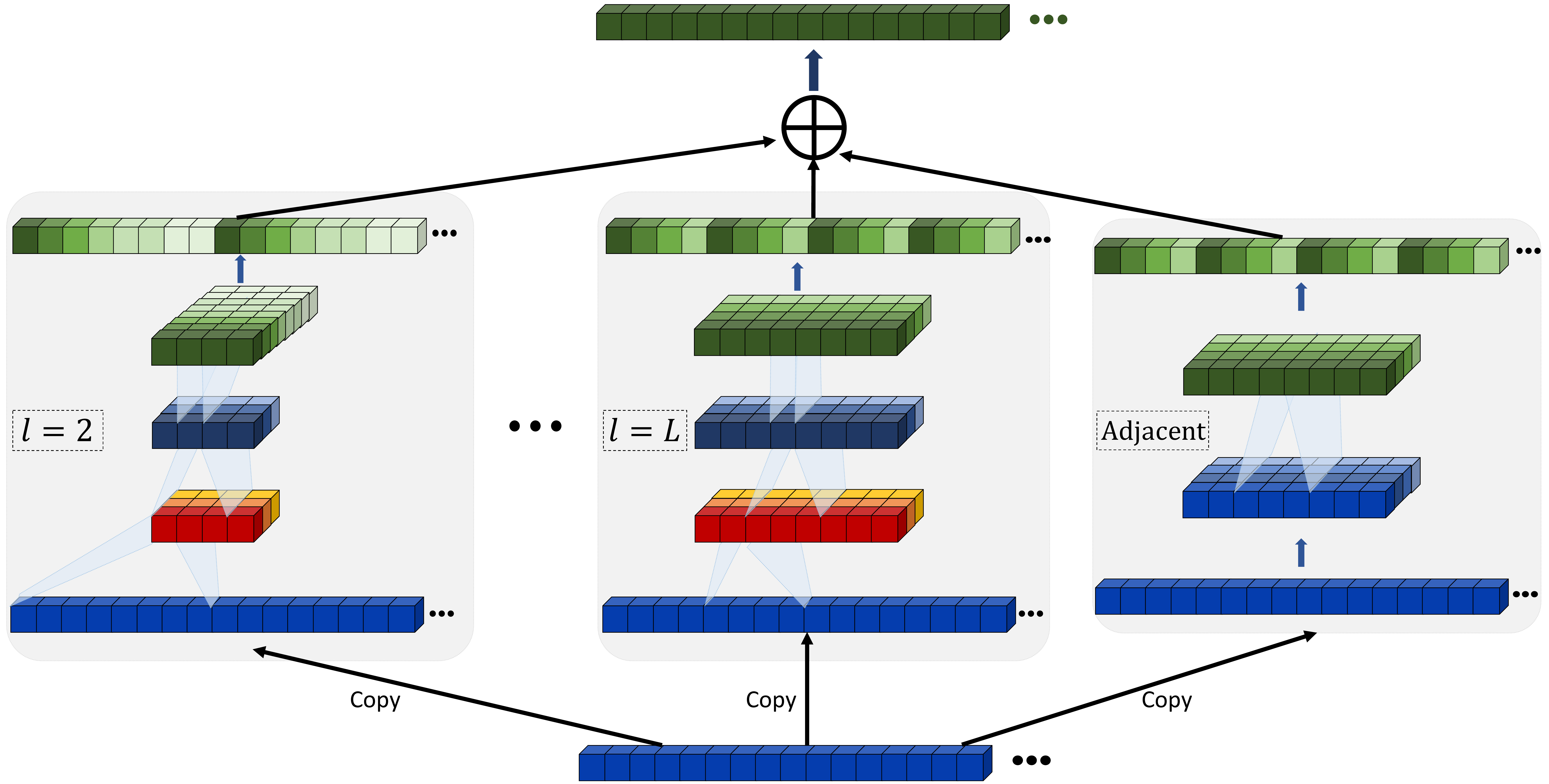}
        \put(17.3, 13.7){\tiny$\LCR$-$\linear$}
        \put(17.3, 19.7){\tiny$\LCK$-$\linear$}
        \put(17.3, 25.7){\tiny$\LCI$-$\linear$}
        \put(57.3, 13.7){\tiny$\LCR$-$\linear$}
        \put(57.3, 19.7){\tiny$\LCK$-$\linear$}
        \put(57.3, 25.7){\tiny$\LCI$-$\linear$}
        \put(89.3, 22.7){\tiny$\LCK$-$\linear$}
        \put(84.8, 30.7){\tiny $\Flatten$}
        \put(84.8, 14.4){\tiny $\Reshape$}
        \put(52.5, 32.3){\tiny $\Flatten$}
        \put(16.5, 32.5){\tiny $\Flatten$}
        \put(53.2, 44.8){\small$\mathsf{sum}$}
    \end{overpic}
    \caption{\label{fig:hnnlinear} Neural network architecture for $\cH$-matrices.}
\end{figure}

Combining \eqref{eq:lcR}, \eqref{eq:lcM} and \eqref{eq:lcAd}, we obtain 
\cref{alg:Hmatrix}, whose architecture is illustrated in \cref{fig:hnnlinear}. 
In particular, \cref{fig:hnnlinear} is the translation to the neural network framework of
\eqref{eq:H-Av} (see \cref{fig:Av}) using the building blocks depicted in \cref{fig:LC}. 

\begin{algorithm}[htb]
\begin{small}
\begin{center}
\begin{minipage}{0.45\textwidth}
\begin{algorithmic}[1]
    \STATE u = 0;
    \FOR {$\ell = 2$ to $L$}
    \STATE $\xi = \LCR[\linear; N,2^{\ell},r](v)$;
    \STATE $\zeta = \LCK[\linear; 2^{\ell},r, r, 2n_b\sps{\ell}+1](\xi)$;
    \STATE $u = u + \LCI[\linear; 2^{\ell}, r, \frac{N}{2^{\ell}}](\zeta)$;
    \ENDFOR
\end{algorithmic}
\end{minipage}
\begin{minipage}{0.5\textwidth}
\begin{algorithmic}[1]
    \setcounter{ALC@line}{6}
    \STATE $\tilde{v}=\Reshape[m,2^L](v)$;
    \STATE $\tilde{u}^{(\ad)}=\LCK\left[\linear; 2^L,m,m,2n^{(\ad)}_{b}+1\right](\tilde{v})$;
    \STATE $u^{(\ad)}=\Flatten(\tilde{u}^{(\ad)})$;
    \STATE $u = u + u^{(\ad)}$;
\end{algorithmic}
\vspace{3mm}
\end{minipage}
\end{center}
\end{small}
\caption{Application of the NN representation of an $\cH$-matrix to a vector $v\in\bbR^N$.}\label{alg:Hmatrix}\end{algorithm}

Moreover, the memory footprints of the neural network architecture and $\mathcal{H}$-matrices are
asymptotically the same with respect the spatial dimension of $u$ and $v$.
This can be readily shown by computing the total number of parameters. For the sake of simplicity, 
we only count the parameters in the weights, ignoring those in the biases. A
direct calculation yields the number of parameters in $\LCR$, $\LCK$ and $\LCI$:
\begin{equation}
    N_p^{\LCR} = N_x\alpha',\quad
    N_p^{\LCK} = N_x\alpha\alpha'w,\quad
    N_p^{\LCI} = N_x\alpha\alpha',
\end{equation}
respectively. Hence, the number of parameters in \cref{alg:Hmatrix} is
\begin{equation}
    \begin{aligned}
    N_p^{\cH} &= \sum_{\ell=2}^L\left( Nr + 2^{\ell}r^2(2n_b\sps{\ell}+1) + 2^{\ell}r\frac{N}{2^{\ell}}\right)
    + 2^Lm^2(2n^{(\ad)}_{b}+1)\\
    &\leq 2LNr+2^{L+1}r^2(2\max_{\ell=2}^L n_b\sps{\ell}+1)+Nm(2n^{(\ad)}_{b}+1)\\
    &\leq 2N\log(N)r+3Nm(2n_b+1) \sim O(N\log(N)),
    \end{aligned}
\end{equation}
where $n_b = \max(n^{(\ad)}_{b}, n_{b}\sps{\ell}, \ell=2,\cdots,L)$, and $r\leq m$ is used.

\subsection{Multi-dimensional case}\label{sec:nD}
Following the previous section, the extension of \cref{alg:Hmatrix} to the $d$-dimensional 
case can be easily deduced using the tensor product of one-dimensional cases. 
Consider $d=2$ below for instance, and the generalization to higher dimensional case will be
straight-forward.
Suppose that we have an IE in 2D given by 
\begin{equation}\label{eq:integral2D}
    u(x) = \int_{\Omega}g(x,y)v(y) \dd y, \quad \Omega=[0,1)\times[0,1) ,
\end{equation}
we discretize the domain
$\Omega$ with a uniform grid with $n=N^2$ ($N=2^Lm$)
discretization points, and let $A$ be the resulting matrix obtained from discretizing \eqref{eq:integral2D}. 
We denote the set of all grid points as
\begin{equation}
    \cI\sps{d,0}=\{(k_1/N, k_2/N): k_1,k_2=0,\dots,N-1\}.
\end{equation}
Clearly $\cI\sps{d,0} = \cI\sps{0}\otimes\cI\sps{0}$, where $\cI\sps{0}$ is defined in \eqref{eq:cI},
and $\otimes$ is tensor product.
At each level $\ell$ ($0\leq \ell\leq L$), we decompose the grid in $4^l$ disjoint boxes as 
$\cI_i\sps{d,\ell}=\cI_{i_1}\sps{\ell}\otimes\cI_{i_2}\sps{\ell}$ for $i_1,i_2=1,\dots,2^l$.
The definition of the lists $\NL$ and $\IL$ can be easily extended. 
For each box $\cI$, $\NL(\cI)$ contains $3$ boxes for 1D case, $3^2$ boxes for 2D case. Similarly, 
the decomposition \eqref{eq:decompose} on the matrix $A$ can easily to extended for this case.
Following the structure of $\cH$-matrices, the off-diagonal blocks of $A\sps{\ell}$ can be approximated as
\begin{equation}
  A\sps{\ell}_{\cI, \cJ}\approx U\sps{\ell}_{\cI}M\sps{\ell}_{\cI,\cJ}(V\sps{\ell}_{\cJ})^{T},
    \quad \cI, \cJ\in\cI\sps{\ell},
    \quad U\sps{\ell}_{\cI}, V\sps{\ell}_{\cJ}\in\bbR^{(N/ 2^{\ell})^2\times r},
    M\sps{\ell}_{\cI,\cJ}\in\bbR^{r\times r}.
\end{equation}

As mentioned before, we can describe the network using tensor or vectors. In what follows we will 
switch between representations in order to illustrate the concepts in a compact fashion. We denote 
an entry of a tensor $T$ by $T_{i,j}$, where $i$ is $2$-dimensional index $i=(i_1,i_2)$.
Using the tensor notations, $U\sps{\ell}$, $V\sps{\ell}$ in \eqref{eq:H-Av}, can be treated as 4-tensors
of dimension $N\times N \times 2^\ell r\times 2^\ell$.
We generalize the notion of band matrix to {\em band tensors}. A band tensor $T$ satisfies that 
\begin{equation}
    T_{i,j}=0,\quad \text{if } |i_1-j_1|>n_{b,1} ~\text{ or }~ |i_2-j_2| > n_{b,2},
\end{equation}
where $n_b=(n_{b,1}, n_{b,2})$ is called the band size for tensor. Thus
\cref{pro:A} can be generalized to tensors yielding the following properties. 
\begin{property}
    The 4-tensors
    \begin{enumerate}
        \item $U\sps{\ell}$ and $V\sps{\ell}$, $\ell=2,\cdots,L$ are block diagonal tensors with
            block size $N / 2^{\ell}\times N / 2^{\ell} \times r\times 1$;
        \item $A^{(\ad)}$ is a block band cyclic tensor with block size $m\times m\times m \times m$
            and band size $n^{(\ad)}_{b}=(1,1)$;
        \item $M\sps{\ell}$, $\ell=2,\cdots,L$ are block band cyclic tensor with block size $r\times
            1\times r\times 1$ and band size $n_{b}\sps{\ell}$, which is $(2,2)$ for $\ell=2$ and $(3,3)$ for
            $\ell\geq3$.
    \end{enumerate}
\end{property}
 
Next, we characterize LC networks for the 2D case. An NN layer for 2D can be represented by a 3-tensor of size $\alpha\times N_{x,1}\times N_{x,2}$, in which $\alpha$ is
the channel dimension and $N_{x,1}$, $N_{x,2}$ are the spatial dimensions. If a layer $\xi$ with
size $\alpha\times N_{x,1}\times N_{x,2}$ is connected to a locally connected layer $\zeta$ with
size $\alpha'\times N_{x,1}'\times N_{x,2}'$, then
\begin{footnotesize}
\begin{equation}\label{eq:lc2d}
    \zeta_{c',i} = \phi\left( \sum_{j=(i-1)s+1}^{(i-1)s+w}\sum_{c=1}^{\alpha}
    W_{c',c;i,j}\xi_{c,j} + b_{c',i}\right),
    \quad i_1=1,\dots,N_{x,1}', i_2=1,\dots,N_{x,2}',~ c'=1,\dots,\alpha',
\end{equation}
\end{footnotesize}where $(i-1)s= ( (i_1-1)s_1, (i_2-1)s_2)$.
As in the 1D case, the channel dimension corresponds to the rank $r$, and the spatial dimensions
correspond to the grid points of the discretized domain. 
Analogously to the 1D case, we define the LC networks 
$\LCR$, $\LCK$ and $\LCI$ and use them to express the four operations in \eqref{eq:hm} which 
constitute the building blocks of the neural network. 
The extension is trivial, the parameters $N_x$, $s$ and $w$ in the one-dimensional LC networks are
replaced by their 2-dimensional counterpart $N_x=(N_{x,1}, N_{x,2})$, $s=(s_1,s_2)$ and
$w=(w_{1},w_{2})$, respectively. 
We point out that $s=w=\frac{N_x}{N_x'}$ for the 1D case is replaced by
$s_j=w_{j}=\frac{N_{x,j}}{N_{x,j}'}$, $j=1,2$ for the 2D case in the definition of LC.

\begin{figure}[ht]
    \centering
    \includegraphics[width=0.3\textwidth]{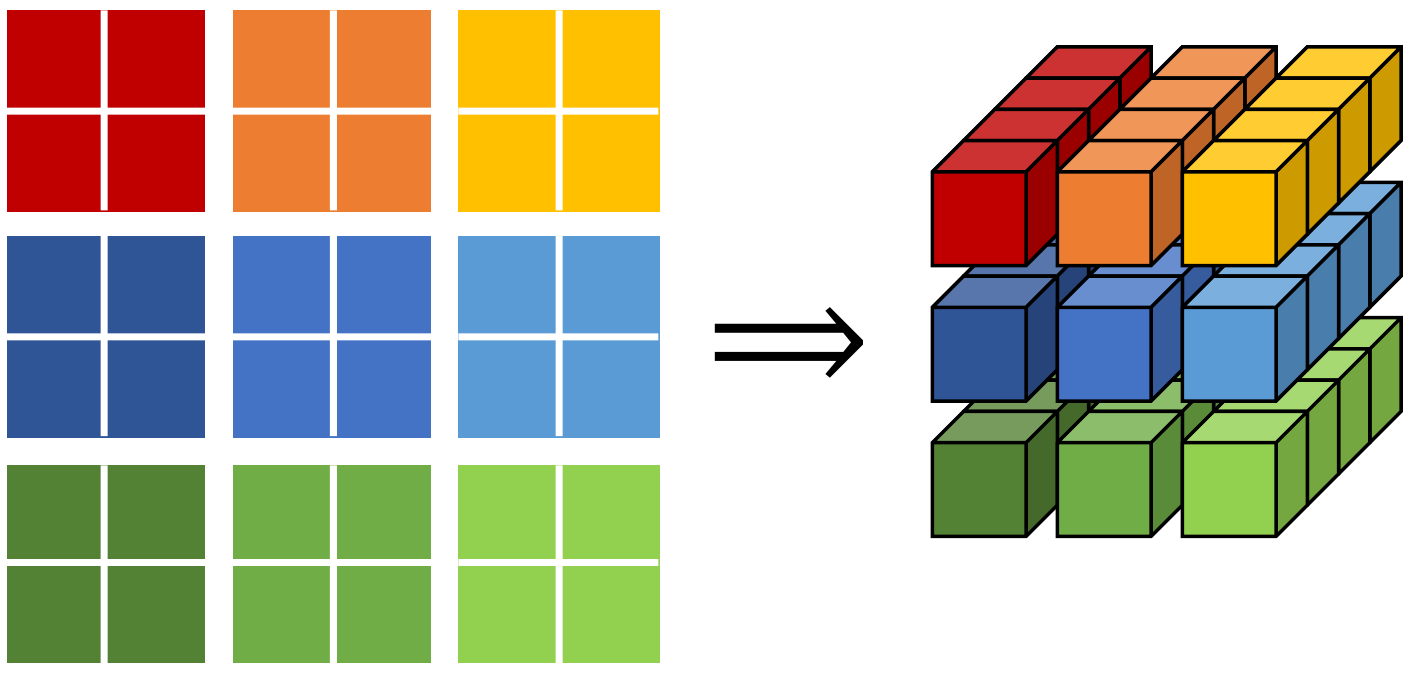}
    \caption{\label{fig:reshape2d}Diagram of $\Reshape[2^2,3,3]$ in \cref{alg:Hmatrix2d}.}
\end{figure}

Using the notations above we extend \cref{alg:Hmatrix} to the 2D case in \cref{alg:Hmatrix2d}. We
crucially remark that the $\Reshape[r^2, n_1, n_2]$ function in \cref{alg:Hmatrix2d} 
is not the usual major column based reshaping. It reshapes a 2-tensor $T$ with size $rn_1\times rn_2$  
to a 3-tensor $S$ with size $r^2\times n_1\times n_2$, by treating the former as a block tensor
with block size $r\times r$, and reshaping each block as a vector following the formula 
$S(k, i,j)=T((i-1)r+k_1, (j-1)r+k_2)$ with $k=(k_1-1)r+k_2$, for 
$k_1,k_2=1,\dots,r$, $i=1,\dots,n_1$ and $j=1,\dots,n_2$. \cref{fig:reshape2d} 
provides an example for the case $\Reshape[2^2,3,3]$. The $\Flatten$ is its inverse. 

\begin{algorithm}[htb]
\begin{small}
\begin{center}
\begin{minipage}{0.9\textwidth}
\begin{algorithmic}[1]
    \STATE u = 0;
    \FOR {$\ell = 2$ to $L$}
    \STATE $\xi = \LCR[\linear; (N, N), (2^{\ell},2^{\ell}),r](v)$;
    \STATE $\zeta = \LCK[\linear; (2^{\ell},2^{\ell}),r, r, (2n_{b,1}\sps{\ell}+1,2n_{b,2}\sps{\ell}+1)](\xi)$;
    \STATE $u = u + \LCI[\linear; (2^{\ell},2^{\ell}), r, \left(\frac{N}{2^{\ell}}\right)^2](\zeta)$;
    \ENDFOR
    \STATE $\tilde{v}=\Reshape[m^2,2^L, 2^L](v)$;
    \STATE $\tilde{u}\sps{\ad}=\LCK\left[\linear; (2^L,2^L),m^2,m^2,(2n\sps{\ad}_{b,1}+1,2n\sps{\ad}_{b,2}+1)\right](\tilde{v})$;
    \STATE $u\sps{\ad}=\Flatten(\tilde{u}\sps{\ad})$;
    \STATE $u = u + u\sps{\ad}$;
\end{algorithmic}
\end{minipage}
\end{center}
\end{small}
\caption{Application of NN architecture for $\cH$-matrices on a vector $v\in\bbR^{N^2}$.}\label{alg:Hmatrix2d}\end{algorithm}

\section{Multiscale neural network}\label{sec:hnn}

In this section, we extend the aforementioned NN architecture to represent a 
\revised{general nonlinear mapping of the form}
{nonlinear generalization of pseudo-differential operators of the form}
\begin{equation}\label{eq:nonlinearmap}
    u = \cM(v),\quad u,v\in\bbR^{N^d}.
\end{equation}
Due to its multiscale structure, we refer to the resulting NN architecture
as the {\em multiscale neural network} (MNN). We consider the one-dimensional case below for simplicity,
and the generalization to higher dimensions follows directly as in \cref{sec:nD}.

\subsection{Algorithm and architecture}\label{sec:mnnAlgorithm}
\begin{algorithm}[ht]
\begin{small}
\begin{center}
\begin{minipage}{0.45\textwidth}
\begin{algorithmic}[1]
    \STATE u = 0;
    \FOR {$\ell = 2$ to $L$}
    \STATE $\xi_0 = \LCR[\linear; N,2^{\ell},r](v)$;
    \FOR {$k=1$ to $K$}
    \STATE $\xi_k= \LCK[\phi; 2^{\ell},r, r,2n_b\sps{\ell}+1](\xi_{k-1})$;
    \ENDFOR
    \STATE $u = u + \LCI[\linear; 2^{\ell}, r, \frac{N}{2^{\ell}}](\xi_{K})$;
    \ENDFOR
\end{algorithmic}
\end{minipage}
\begin{minipage}{0.5\textwidth}
\begin{algorithmic}[1]
    \setcounter{ALC@line}{8}
    \STATE $\xi_0=\Reshape[m,2^L](v)$;
    \FOR {$k=1$ to $K-1$}
    \STATE $\xi_k=\LCK[\phi; 2^L, m, m,2n^{(\ad)}_{b}+1](\xi_{k-1})$;    \ENDFOR
    \STATE \add{$\xi_K=\LCK[\linear; 2^L, m, m,2n^{(\ad)}_{b}+1](\xi_{K-1})$;}\label{alg:LCKad}
        \STATE $u^{(\ad)}=\Flatten(\xi_K)$;
    \STATE $u = u + u^{(\ad)}$;
\end{algorithmic}
\end{minipage}
\end{center}
\end{small}
\caption{Application of multiscale neural network to a vector $v\in\bbR^N$.}
\label{alg:mnn}
\end{algorithm}

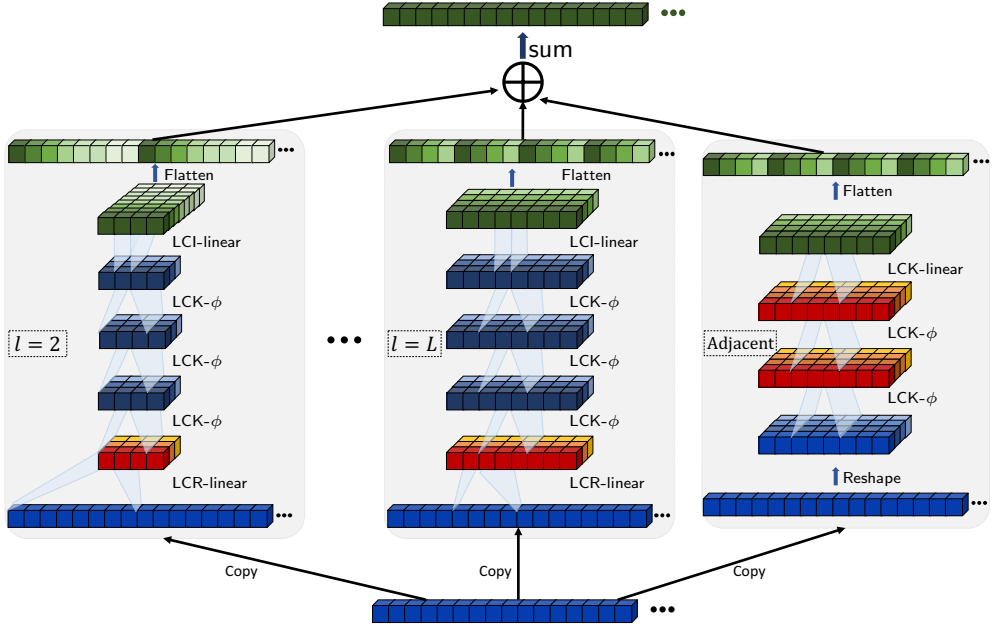
\begin{figure}[htb]
    \centering
    \begin{overpic}[width=0.8\textwidth]{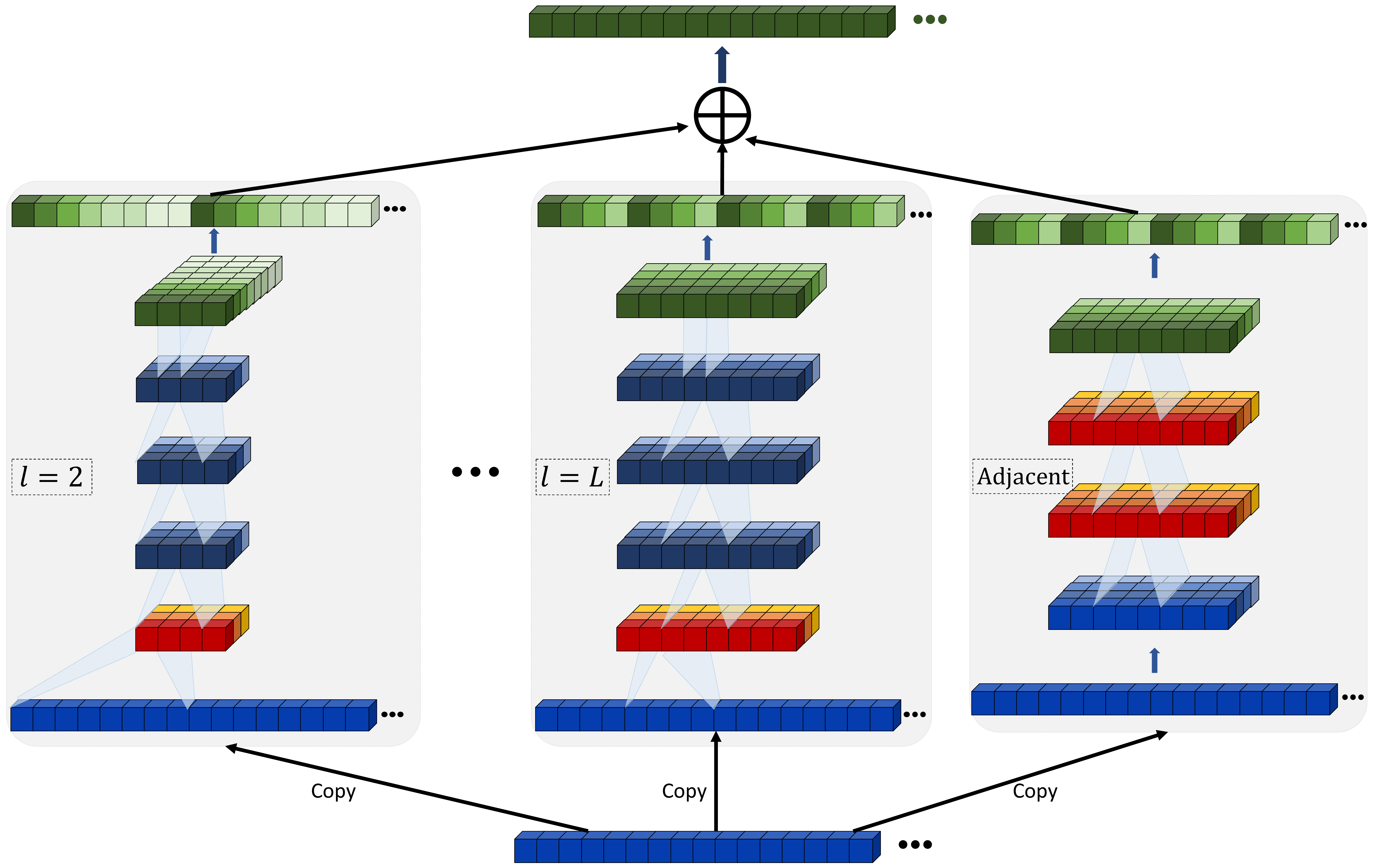}
        \put(17.3, 13.9){\tiny$\LCR$-$\linear$}
        \put(17.3, 20.2){\tiny$\LCK$-$\phi$}
        \put(17.3, 26.2){\tiny$\LCK$-$\phi$}
        \put(17.3, 32.2){\tiny$\LCK$-$\phi$}
        \put(17.3, 38.2){\tiny$\LCI$-$\linear$}
        \put(16.5, 45.0){\tiny $\Flatten$}
        \put(57.3, 13.9){\tiny$\LCR$-$\linear$}
        \put(57.3, 20.2){\tiny$\LCK$-$\phi$}
        \put(57.3, 26.2){\tiny$\LCK$-$\phi$}
        \put(57.3, 32.2){\tiny$\LCK$-$\phi$}
        \put(57.3, 38.2){\tiny$\LCI$-$\linear$}
        \put(56.5, 45.0){\tiny $\Flatten$}
        \put(89.3, 22.5){\tiny$\LCK$-$\phi$}
        \put(89.3, 29.0){\tiny$\LCK$-$\phi$}
        \put(89.3, 35.5){\tiny$\LCK$-$\linear$}
        \put(84.8, 14.4){\tiny $\Reshape$}
        \put(84.8, 43.3){\tiny $\Flatten$}
        \put(53.2, 57.4){\small$\mathsf{sum}$}
    \end{overpic}
    \caption{\label{fig:mnn} Multiscale neural network architecture for nonlinear mappings, which is an
    extension of the neural network architecture for $\cH$-matrices \cref{fig:hnnlinear}.  
    $\phi$ is an activation function.}
\end{figure}

NN can represent nonlinearities by choosing the activation function, $\phi$, to be nonlinear, such as
ReLU or sigmoid.  The range of the activation function also imposes constraints on the output of the
NN. For example, the range of ``ReLU'' in $[0,\infty)$ and the range of the sigmoid function is
$[0,1]$. Thus,
the last layer is often chosen to be a linear layer to relax such constraint. 
\cref{alg:Hmatrix} is then revised to \cref{alg:mnn}, and the architecture is
illustrated in \cref{fig:mnn}.
We remark that the nonlinear activation function is only used in the $\LCK$ network.  
The $\LCR$ and $\LCI$ networks in \cref{alg:Hmatrix} are still treated as restriction and
interpolation operations between coarse grid and fine grid, respectively, so we use the linear
activation functions in these layers. Particularly, we also use the linear activation function for
the last layer of the adjacent part, 
which is marked in line \ref{alg:LCKad} in \cref{alg:mnn}.

As in the linear case, we calculate the number of parameters of MNN and obtain (neglecting the
number of parameters in $b$ in \eqref{eq:lc})
\begin{equation}
    \begin{aligned}
        N_p^{\mathrm{MNN}} &= \sum_{\ell=2}^L\left( Nr +
        K2^{\ell}r^2(2n_b\sps{\ell}+1) + 2^{\ell}r\frac{N}{2^{\ell}}\right)
        + K2^Lm^2(2n^{(\ad)}_{b}+1)\\
        &\leq 2LNr+K2^{L+1}r^2(2\max_{\ell=2}^L n_b\sps{\ell}+1)+NKm(2n^{(\ad)}_{b}+1)\\
        &\leq 2N\log(N)r+3NKm(2n_b+1).
    \end{aligned}
\end{equation}

\subsection{Translation-invariant case}\label{sec:equivariant}
For the linear case \eqref{eq:integral}, if the kernel is {\em translation invariant}, \ie $g(x,y) =
g(x-y)$, then the matrix $A$ is a Toeplitz matrix. Then the matrices $M\sps{\ell}$ and
$A^{(\ad)}$ are Toeplitz matrices and all matrix blocks of $U\sps{\ell}$ (resp. $V\sps{\ell}$) can
be represented by the same matrix. 
This leads to the 
{\em convolutional neural network} (CNN) as
\begin{equation}\label{eq:cnn}
    \zeta_{c',i} = \phi\left( \sum_{j=(i-1)s+1}^{(i-1)s+w}\sum_{c=1}^{\alpha}
    W_{c',c;j}\xi_{c,j} + b_{c'}\right),
    \quad i=1,\dots,N_{x}', ~ c'=1,\dots,\alpha'.
\end{equation}
Compared to the LC network, the only difference is that the parameters $W$ and
$b$ are independent of $i$. Hence, inheriting the
definition of $\LCR$, $\LCK$ and $\LCI$, we define the layers $\CNNR$, $\CNNK$ and $\CNNI$,
respectively. By replacing the LC layers in \cref{alg:Hmatrix} by the corresponding CNN
layers, we obtain the neural network architecture for the translation invariant kernel.

For the nonlinear case, the translation invariant kernel for the linear case can be extended to
kernels that are {\em equivariant to translation}, \ie for any translation $\cT$,
\begin{equation}\label{eq:equivariant}
    \cT\cM(v) = \cM(\cT v).
\end{equation}
For this case, all the LC layers in \cref{alg:mnn} can be replaced by its corresponding CNN layers.
The number of parameters of $\CNNR$, $\CNNK$ and $\CNNI$ are
\begin{equation}
    N_p^{\CNNR} = \frac{N_x}{N_x'}\alpha',\quad
    N_p^{\CNNK} = \alpha\alpha'w,\quad
    N_p^{\CNNI} = \alpha\alpha'.
\end{equation}
Thus, the number of parameters in \cref{alg:mnn} using CNN is
\begin{equation}
    \begin{aligned}
        N_{p,CNN}^{\mathrm{MNN}} &= \sum_{\ell=2}^L\left(
        r\frac{N}{2^{\ell}} + Kr^2(2n_b\sps{\ell}+1) + r\frac{N}{2^{\ell}}\right)
        + Km^2(2n^{(\ad)}_{b}+1)\\
                &\leq rN + (r^2\log(N)+m^2)(2n_b+1)K.
    \end{aligned}
\end{equation}

 \section{Numerical results}\label{sec:application}
In this section we discuss the implementation details of MNN. We demonstrate the
accuracy of the MNN architecture using two nonlinear problems: the nonlinear Schr{\"o}dinger equation
(NLSE), and the Kohn-Sham map (KS map) in the Kohn-Sham density functional theory (KSDFT).

\subsection{Implementation}
Our implementation of MNN uses Keras \cite{keras}, a high-level application programming interface (API) running,
in this case, on top of TensorFlow \cite{tensorflow} (a library of tools for training neural networks).
The loss function is chosen as the mean squared relative error, in which the relative error is defined
with respect to $\ell^2$ norm as
\begin{equation}\label{eq:relativeerror}
  \epsilon = \frac{||u-u_{NN}||_{\ell^2}}{||u||_{\ell^2}},
\end{equation}
where $u$ is the target solution generated by a numerical discretization of the PDEs and $u_{NN}$ is the
predicted solution by MNN.  The optimization is performed using the NAdam optimizer
\cite{dozat2015incorporating}. The weights and biases in MNN are initialized randomly from the normal
distribution and the batch size is always set between $1 / 100$th and $1 / 50$th of the number of train samples.

In all the tests, the band size is chosen as $n_{b,\ad}=1$ and $n_{b}\sps{l}$ is $2$ for $l=2$ and
$3$ otherwise. The nonlinear activation function is chosen as ReLU.  All the test are run on GPU
with data type \verb|float32|. All the numerical results are the best results by repeating the
training a few times, using different random seeds.
The selection of parameters $r$ (number of channels), $L$ ($N=2^Lm$) and $K$ (number of layers in
\cref{alg:mnn}) is problem dependent.

\subsection{NLSE with inhomogeneous background potential}\label{sec:nlse}
The nonlinear Schr{\"o}dinger equation (NLSE) is widely used in quantum physics to describe the
single particle properties of the Bose-Einstein condensation phenomenon
\cite{pitaevskii1961vortex,anglin2002bose}.
Here we study the NLSE with inhomogeneous background potential $V(x)$:
\begin{equation}\label{eq:nlse}
    \begin{aligned}
        &-\Delta u(x) + V(x)u(x) + \beta u(x)^3=E u(x),\quad x\in [0,1)^d,\\
        &\text{ s.t.} \int_{[0,1]^d}u(x)^2\dd x=1, \text{ and } \int_{[0,1]^d}u(x)\dd x>0,
    \end{aligned}
\end{equation}
with period boundary condition. We aim to find its ground state denoted by $u_G(x)$. We take a strongly nonlinear case
$\beta=10$ in this work and thus consider a defocusing cubic Schr{\"o}dinger equation. Due to the
cubic term, an iterative method is required to solve \eqref{eq:nlse} numerically. We employ the
method in \cite{bao2004computing} for the numerical solution, which solves a time-dependent NLSE by a
normalized gradient flow. The MNN is used to learn the map from the background potential to the
ground state
\begin{equation}
    V(x) \rightarrow u_G(x).
\end{equation}
This map is equivariant to translation, and thus MNN is implemented using the CNN layers.
\add{The constraints in \eqref{eq:nlse} can be guaranteed by adding a post-correction in the
network as
\begin{equation}
    u = \mathrm{const} \times \frac{\hat{u}}{\|\hat{u}\|_2},
\end{equation}
where $\hat{u}$ is the prediction of the neural network.
}In the following, we study the performance of MNN on 1D and 2D cases.

\subsubsection{One-dimensional case}
For the one-dimensional case, the number of discretization points is $N=320$, and we set $L=6$ and
$m=\frac{N}{2^L}=5$.
\add{The potential $V$ is chosen as
\begin{equation}\label{eq:V}
    V = -20\exp(I_F(v)),
\end{equation}
where $v\in \mathcal{N}(0,1)^{40}$ and $I_F: \bbR^{40}\to\bbR^{320}$ is the Fourier interpolation operator.
}In all the tests, the number of test samples is the same as that the number of
train samples if not properly specified. We perform numerical experiments to study the behavior of
MNN for different number of channels
$r$, different number of $\CNNK$ layers $K$, and different number of training samples \Ntrainsample.
\add{All the networks are trained using 5000 epochs.}

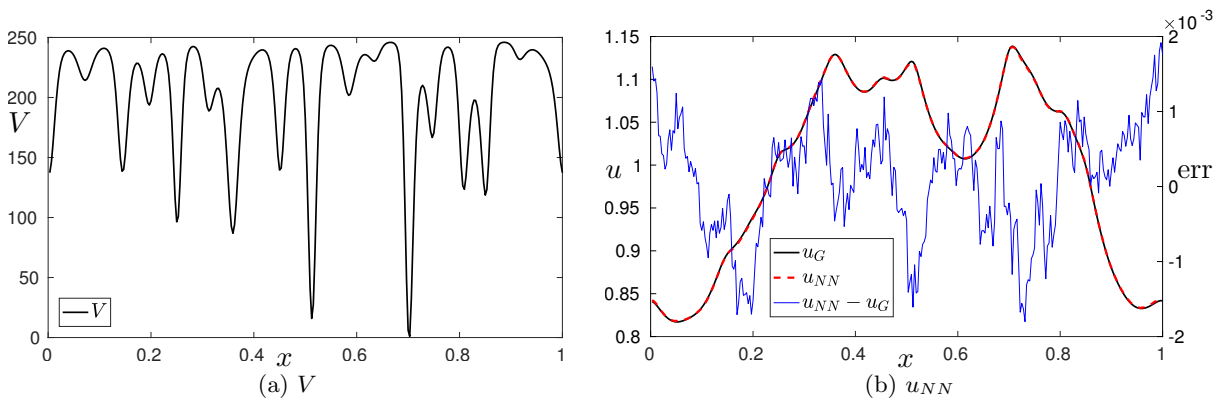
\begin{figure}[ht]
    \centering
    \subfloat[$V$]{
    \begin{overpic}[trim={0mm 0mm 0mm 0mm}, clip, width=0.45\textwidth]{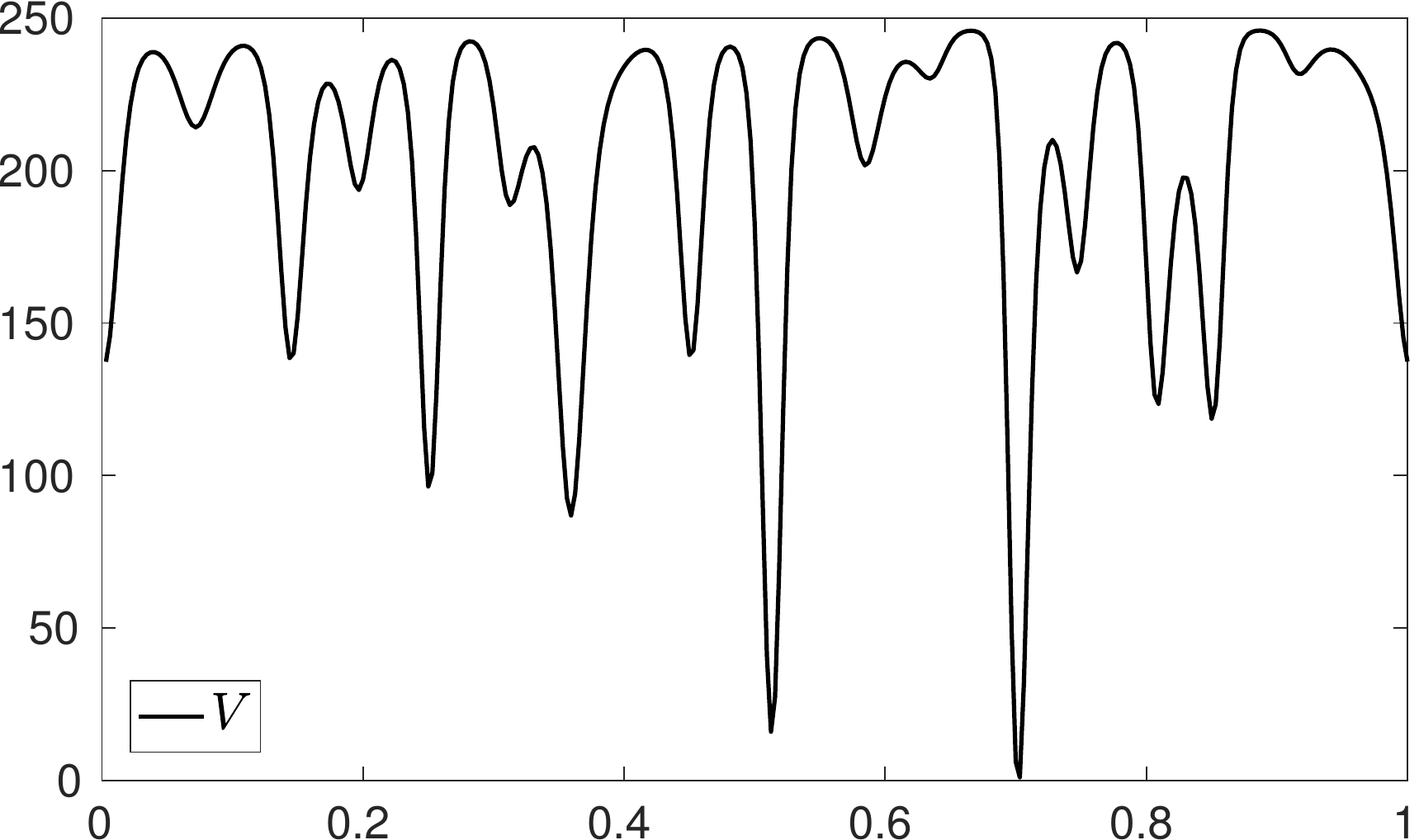}
        \put(48, -1){$x$}
        \put( 0, 41){$V$}
    \end{overpic}
    }~
    \subfloat[$u_{NN}$]{
    \begin{overpic}[trim={0mm 0mm 0mm 0mm}, clip, width=0.49\textwidth]{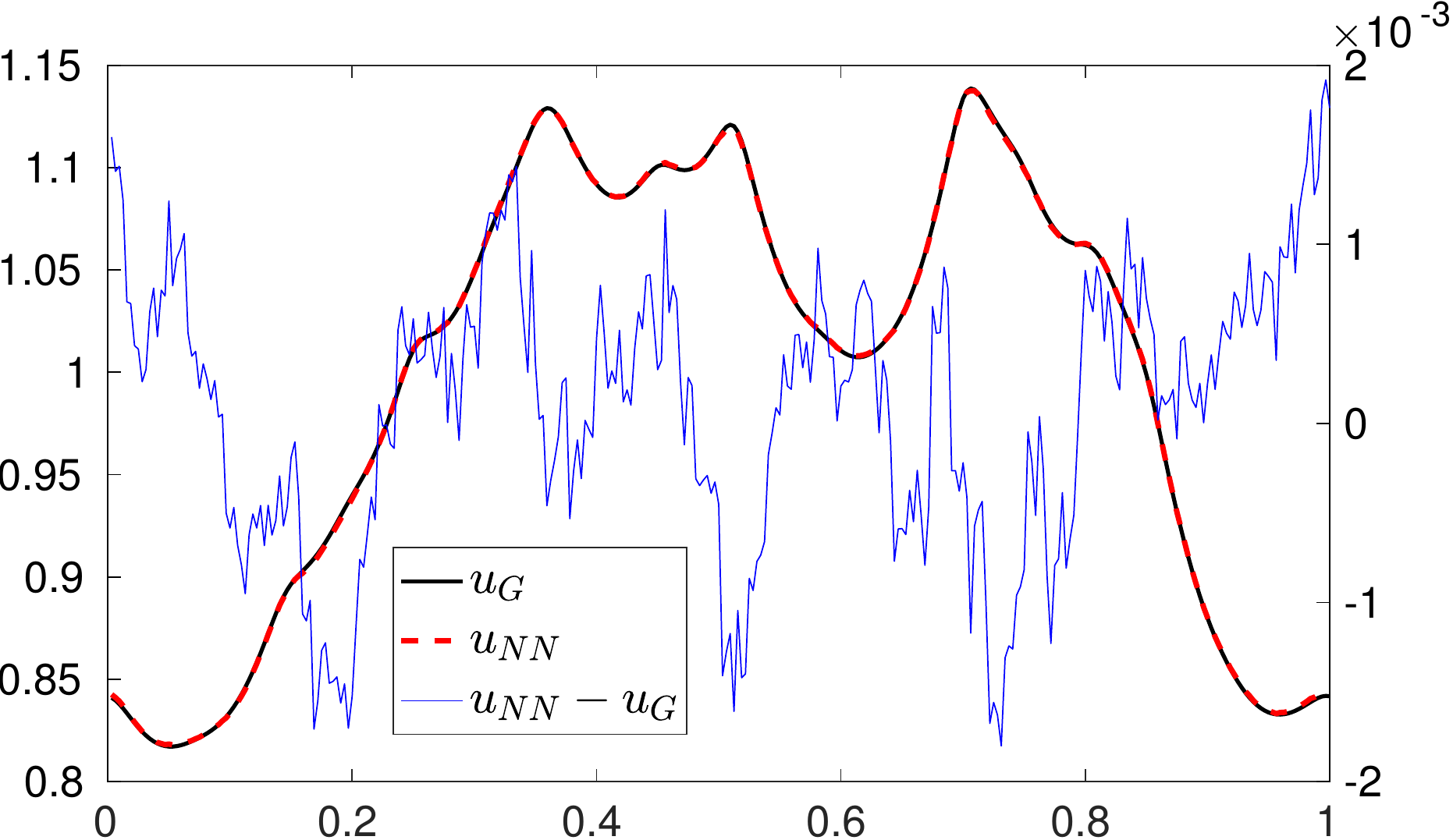}
        \put(48, -1){$x$}
        \put( 0, 30){$u$}
        \put(94, 30){$\mathrm{err}$}
                                    \end{overpic}
    }
    \caption{\label{fig:nlse1d} \add{A sample of the potential $V$ \eqref{eq:V} in the test set and its
    corresponding solution $u_G$, predicted solution $u_{NN}$ and its error with respect to the
    $u_G$ by MNN with $r=8$ and $K=7$ for 1D NLSE.}}
\end{figure}

\begin{table}[htb]
    \centering
    \begin{tabular}{cccc}
      \hline
      \Ntrainsample &	 \Ntestsample &	 Training error 	&	 Validation error\\ \hline\hline
      500	    &	5000	&	8.9e-3	&	1.7e-3\\ \hline
      1000	&	5000	&	9.4e-4	&	1.2e-3\\ \hline
      5000	&	5000	&	8.1e-4	&	8.4e-4\\ \hline
            20000	&	20000	&	8.0e-4	&	8.1e-4\\ \hline
    \hline
    \end{tabular}
    \caption{\label{tab:nlse1d_samples}Relative error in approximating the ground state of NLSE for
    different number of samples \Ntrainsample for 1D case with $r=8$ and $K=7$.}
\end{table}

\begin{table}[ht]
    \centering
    \begin{tabular}{cccc}
      \hline
      $r$	&	 \Nparams    &	 Training error 	&	 Validation error\\ \hline\hline
      2    & 2339    & 2.6e-3    & 2.6e-3\\\hline
      4    & 5811    & 1.2e-3    & 1.2e-3\\\hline
      6    & 11131    & 9.9e-4    & 1.0e-3\\\hline
      8    & 18299    & 8.1e-4    & 8.4e-4\\\hline
                            \hline
    \end{tabular}
    \caption{\label{tab:nlse1d_channel}Relative error in approximating the ground state of NLSE for
    different number of channels $r$ for 1D case with $K=7$ and \Ntrainsample $=5000$.}
\end{table}

\begin{table}[htb]
    \centering
    \begin{tabular}{cccc}
      \hline
      $K$	&	 \Nparams	&	 Training error 	&	 Validation error\\ \hline\hline
      1    & 4907   & 8.2e-3    & 8.2e-3\\\hline
      3    & 9371    & 1.6e-3    & 1.7e-3\\\hline
      5    & 13835    & 1.1e-3    & 1.1e-3\\\hline
      7    & 18299    & 8.1e-4    & 8.4e-4\\\hline
                            \hline
    \end{tabular}
    \caption{\label{tab:nlse1d_layers}Relative error in approximating the ground state of NLSE for
    different number of $\CNNK$ layers $K$ for 1D case with $r=8$ and \Ntrainsample $=5000$.}
\end{table}

\begin{table}
    \centering
    \begin{tabular}{ccccc}
        \hline
                        $\alpha$ & $K$ & $N_{\mathrm{params}}$ & Training error & Validation error \\
        \hline \hline
          8    & 13&     18097  &     4.6e-3&     4.6e-3\\\hline
         10    & 13&     28121  &     3.9e-3&     3.9e-3\\\hline
         12    & 13&     40345  &     3.5e-3&     3.5e-3\\\hline
         14    & 13&     54769  &  3.8e-3  &  3.8e-3  \\\hline
         12    & 15&     47569  &  4.0e-3  &  4.0e-3  \\\hline
        \hline \hline
    \end{tabular}
    \caption{\label{tab:nlse1d_CNN}\add{Relative error in approximating the ground state of NLSE for the
    CNN with different number of channel $\alpha$ and different number of layers $K$ with windows
    size to be $25$ for 1D case with \Ntrainsample$=5000$.}}
\end{table}

Usually, the number of samples should be greater than that of parameters to avoid overfitting.
However, in neural networks it has been consistently found that the number of samples can be less than
that of parameters~\cite{ZhangBengioHardtEtAl2016,zhang2017beyond}
We present the numerical results for different \Ntrainsample with $K=7$ and $r=8$ in
\cref{tab:nlse1d_samples}. In this case, the number of parameters is \Nparams $=18299$.
\add{The validation error is slightly larger than the training error even for the case \Ntrainsample
$=500$ and is approaching to the training error as the \Ntrainsample increases.
For the case \Ntrainsample $=5000$, the validation error is close to the training error,
thus there is no overfitting, and the errors are only slightly larger than that when \Ntrainsample
$=20000$.}
This allows us to train MNN with \Ntrainsample $<$ \Nparams. This feature is particularly useful for
high-dimensional cases, given that in such cases \Nparams is usually very large and
generating large amount of samples can be prohibitively expensive.

\cref{tab:nlse1d_channel} presents the numerical results for different number of channels, $r$, (\ie
the rank of the $\cH$-matrix) with $K=7$ and \Ntrainsample $=5000$. As $r$ increases, we find that
the error first consistently decreases and then stagnates. We
use $r=8$ for the 1D NLSE below to balance between efficiency and accuracy.

Similarly, \cref{tab:nlse1d_layers} presents the numerical results for different number of
$\CNNK$ layers, $K$, with $r=8$ and \Ntrainsample $=5000$. The error consistently decreases with
respect to the increase of $K$, as NN can represent increasingly more complex
functions with respect to the depth of the network. However, after a certain threshold, adding more
layers provides very marginal gains in accuracy. In practice, $K=7$ is a good choice for the NLSE for 1D case.

\add{
    We also compare the MNN with an instance of the classical convolutional neural networks (CNN).
    The CNN architecture used in \cite{fan2018mnnh2} is adopted as the reference architecture.
    \Cref{tab:nlse1d_CNN} presents the setup of the networks and the training and validation errors
    for CNN. Clearly, by comparing the results of MNN in \cref{tab:nlse1d_channel,tab:nlse1d_layers}
    with the results of CNN in \cref{tab:nlse1d_CNN}, one can observe that the MNN not only reduces
    the number of parameters, but also improve the accuracy.
}

throughout the results shown in \cref{tab:nlse1d_channel,tab:nlse1d_layers}, the validation
errors are very close to the corresponding training errors, thus no overfitting is observed.
\cref{fig:nlse1d} presents a sample for the potential $V$ and its corresponding
solution and prediction solution by MNN. We can observe that the prediction solution agrees with
the target solution very well.

\subsubsection{Two-dimensional case}
For the two-dimensional case, the number of discretization points is $n=80\times 80$, and we set $L=4$
and $m=5$.
\add{The potential $V$ is chosen as
\begin{equation}\label{eq:V2d}
    V = -20\exp(I_F(v)),
\end{equation}
where $v\in \mathcal{N}(0,1)^{10\times 10}$ and $I_F: \bbR^{10\times10}\to\bbR^{80\times80}$ is the
two-dimensional Fourier interpolation operator.}In all the tests the number of test data is the same as that of the train date.
We perform several simulations to study the behavior of MNN for different number of channels, $r$,
and different number of $\CNNK$ layers ,$K$.
As discussed in 1D case, MNN allows \Ntrainsample $<$ \Nparams.
\add{In all the tests, the number of samples is always $20000$ and the networks are trained using $5000$
epochs. Numerical test shows no overfitting for all the test.}

\begin{figure}[ht]
    \centering
    \subfloat[$V$]{
    \includegraphics[clip,width=0.35\textwidth]{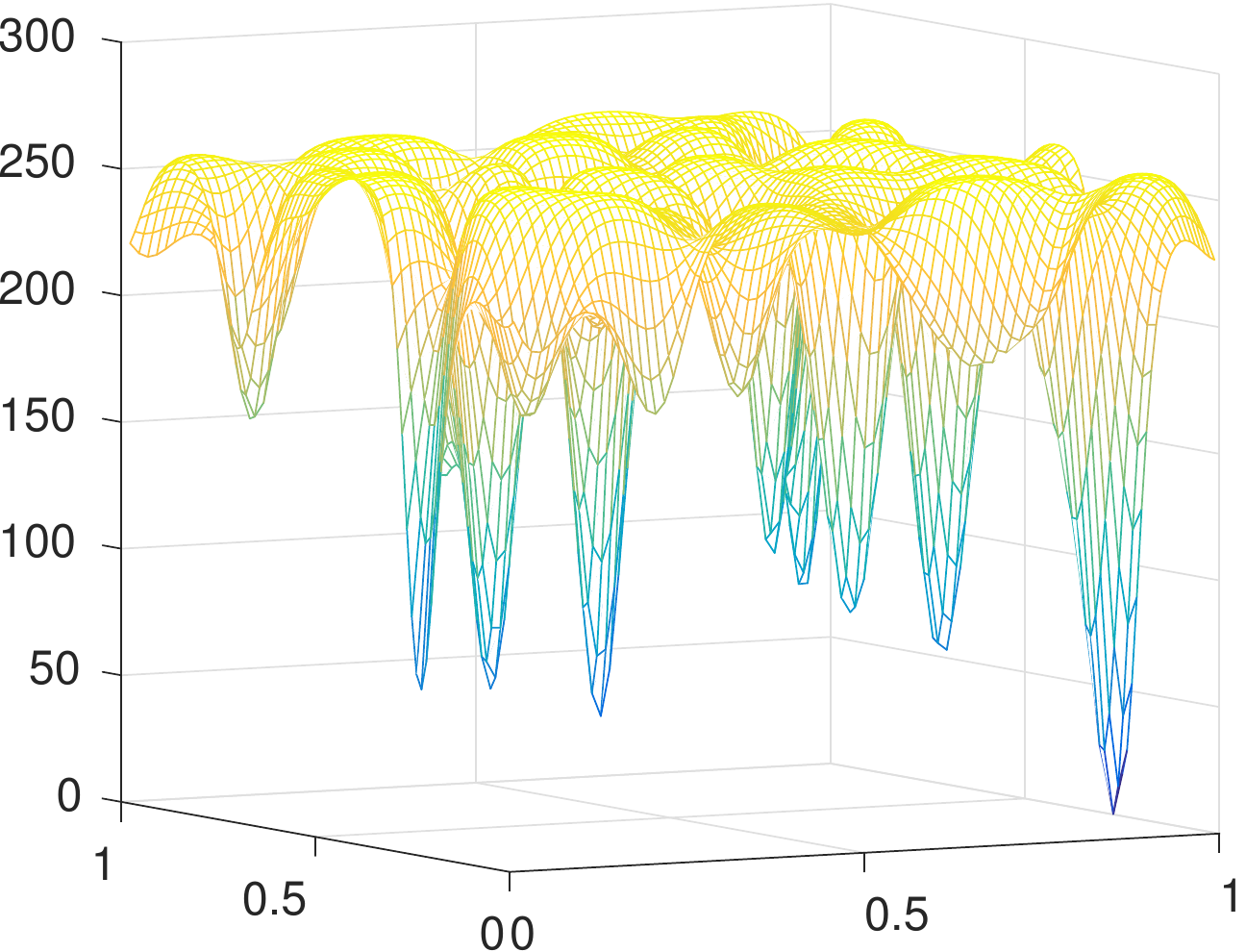}
    }\hspace{0.1\textwidth}
    \subfloat[$u_G$]{
    \includegraphics[clip,width=0.35\textwidth]{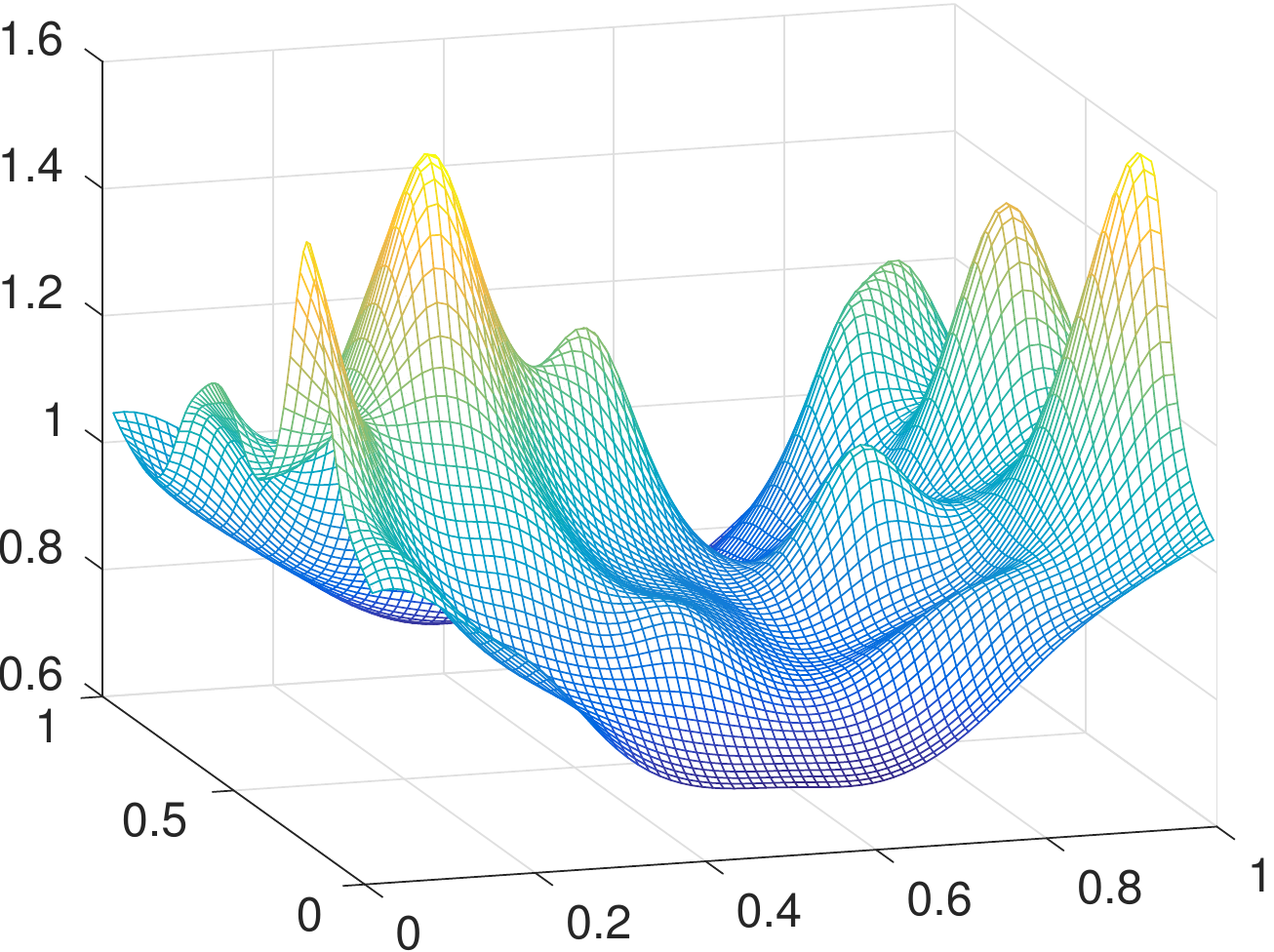}
    }\\
    \subfloat[$u_{NN}$]{
    \includegraphics[clip,width=0.35\textwidth]{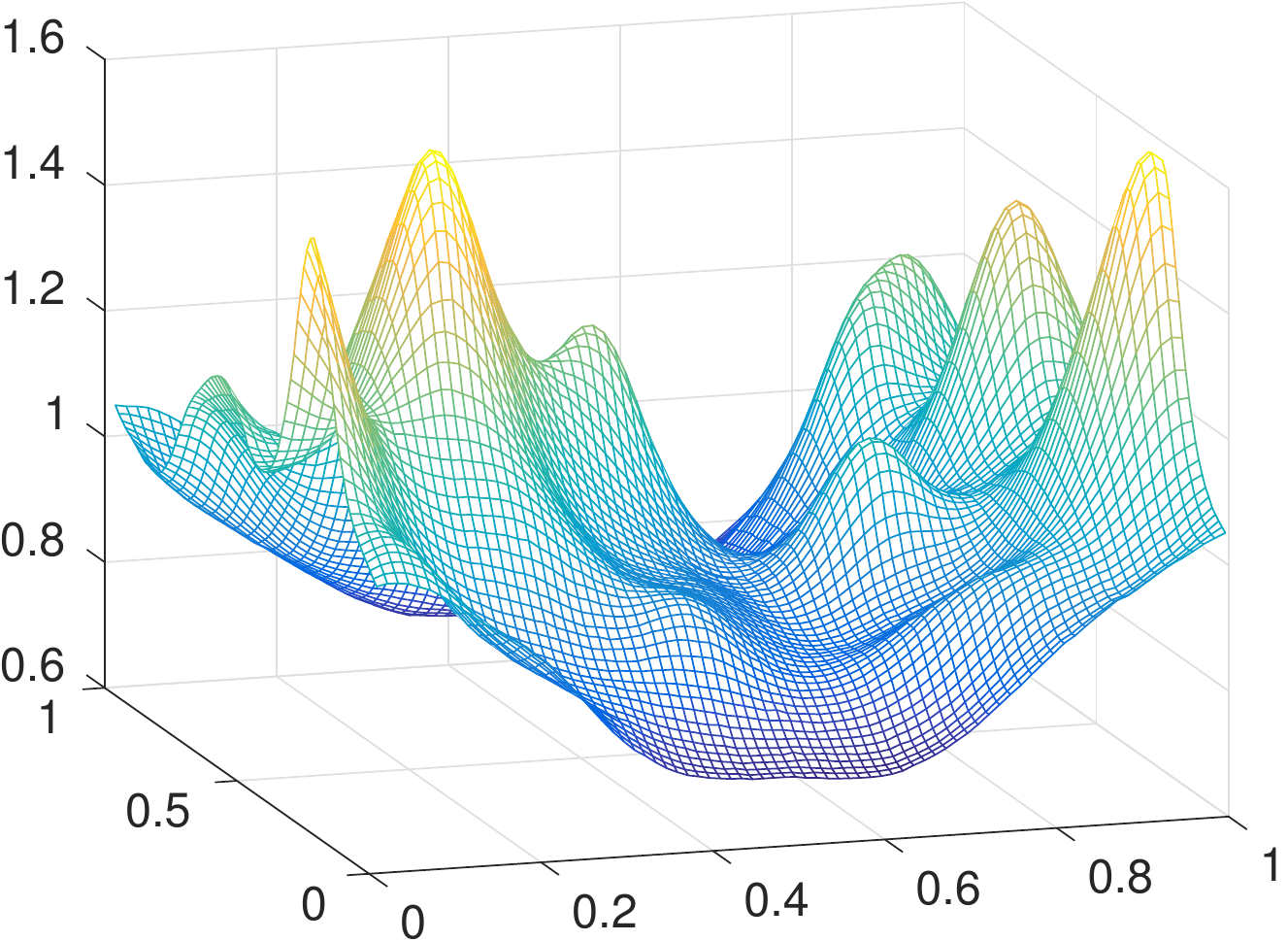}
    }\hspace{0.1\textwidth}
    \subfloat[$u_{NN}-u_G$]{
    \includegraphics[clip,width=0.35\textwidth]{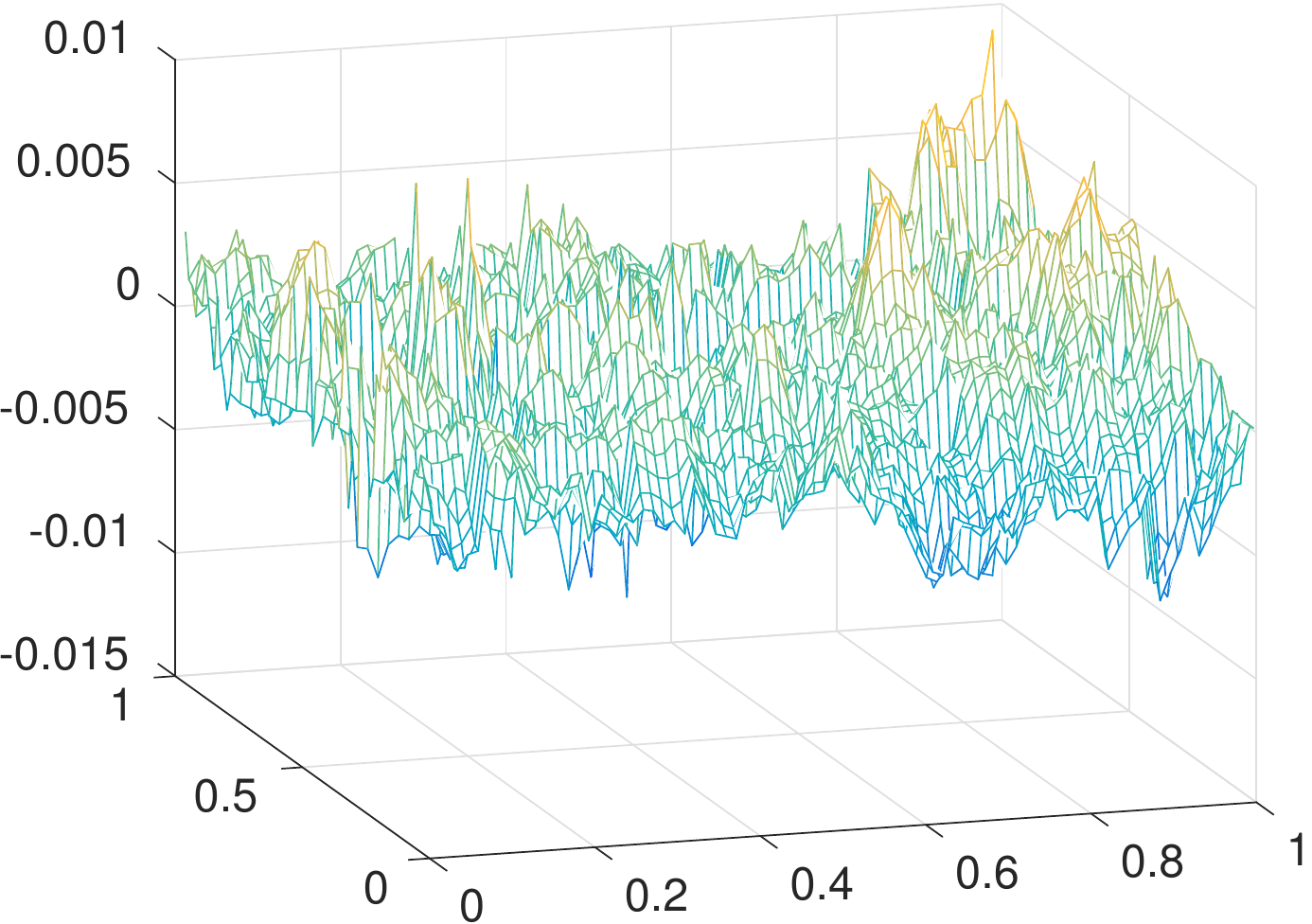}
    }
    \caption{\label{fig:nlse2d}\add{A sample of the potential $V$ \eqref{eq:V2d} in the
    test set and its corresponding solution $u_G$, predicted solution $u_{NN}$ by MNN $K=5$ and $r=6$ and its error
    with respect to $u_G$ for 2D NLSE.}}
\end{figure}

\begin{table}[htb]
  \centering
  \begin{tabular}{cccc}
    \hline
    $r$	&	 \Nparams &	 Training error 	&	 Validation error\\ \hline\hline
    2	&   45667	&   6.1e-3	&   6.1e-3\\ \hline
    6	&   77515	&   2.4e-3	&   2.4e-3\\ \hline
    10	&   136915	&   2.9e-3	&   2.9e-3\\ \hline
                \hline
  \end{tabular}
  \caption{\label{tab:nlse2d_channel}Relative error in approximating the ground state of NLSE for
    different number of channels $r$ for 2D case with $K=7$ and \Ntrainsample $=20000$.}
\end{table}

\begin{table}[htb]
  \centering
  \begin{tabular}{cccc}
    \hline
    $K$	&	 \Nparams &	 Training error 	&	 Validation error\\ \hline\hline
	3	&   37131	&   5.6e-3	&   5.6e-3\\ \hline
	5	&   57323	&   3.8e-3	&   3.8e-3\\ \hline
	7	&   77515	&   2.4e-3	&   2.4e-3\\ \hline
                    \hline
  \end{tabular}
  \caption{\label{tab:nlse2d_layers}Relative error in approximating the ground state of NLSE for
    different number of $\CNNK$ layers $K$ for 2D case with $r=6$ and \Ntrainsample $=20000$.}
\end{table}

\cref{tab:nlse2d_channel,tab:nlse2d_layers} present the numerical results for
different number of channels, $r$, and different number of $\CNNK$ layers ,$K$, respectively.
Similarly to the 1D case, the choice of parameters $r=6$ and $K=7$ also yield accurate
results in the 2D case.
\cref{fig:nlse2d} presents a sample of the potential $V$ in the test set and its corresponding
solution, prediction solution and its error with respect to the reference solution.

\subsection{Kohn-Sham map}\label{sec:KSMap}
Kohn-Sham density functional theory~\cite{HohenbergKohn1964,KohnSham1965} is the most widely used
electronic structure theory. It requires the solution of the following set of nonlinear eigenvalue
equations (real arithmetic assumed for all quantities):
\begin{equation}\label{eqn:KSDFT}
  \begin{split}
    & \left(-\frac{1}{2} \Delta + V[\rho](x)\right) \psi_{i}(x) =
    \varepsilon_{i} \psi_{i}(x), \,\, x \in \Omega=[-1,1)^d \\
    & \int_{\Omega} \psi_{i}(x) \psi_{j}(x) \ud x = \delta_{ij},
    \quad \rho(x) = \sum_{i=1}^{n_e} |\psi_i(x)|^2.
   \end{split}
\end{equation}
Here $n_e$ is the number of electrons (spin degeneracy omitted), $d$ is the spatial dimension,
and $\delta_{ij}$ stands for the
Kronecker delta. In addition, all eigenvalues $\{\varepsilon_{i}\}$ are real and ordered
non-decreasingly, and $\rho(x)$ is the electron density,
which satisfies the constraint
\begin{equation}
  \rho(x)\ge 0,\quad \int_{\Omega} \rho(x) \dd x = n_{e}.
  \label{eqn:rho_constraint}
\end{equation}
The Kohn-Sham equations~\eqref{eqn:KSDFT} need to be solved self-consistently,
which can also viewed
as solving the following fixed point map
\begin{equation}
  \rho = \mathcal{F}_{\text{KS}}[V[\rho]].
  \label{eqn:discrete_KS_map}
\end{equation}
Here the mapping $\mathcal{F}_{\text{KS}}[\cdot]$ from $V$ to $\rho$ is called
the Kohn-Sham map, which for a fixed potential is reduced to a linear eigenvalue
problem, and it constitues the most computationally intensive step for solving
\eqref{eqn:KSDFT}.
We seek to approximate the Kohn-Sham map using a multiscale neural network, whose output was regularized so it satisfies \eqref{eqn:rho_constraint}.

In the following numerical experiments the potential, $V$, is given by
\begin{equation} \label{eqn:gaussian_wells}
    V(x) = -\sum_{i=1}^{n_g} \sum_{j\in \mathbb{Z}^d} c_i \exp\left(-\frac{(x - r_i-2j)^2}{2 \sigma^2}\right),
    \qquad  x \in [-1,1)^d,
\end{equation}
where $d$ is the dimension and $r_i \in [-1,1)^d$. We set \add{$\sigma = 0.05$ or $\sigma = 0.15$ for 1D} 
and $\sigma=0.2$ for 2D. The coefficients $c_i$ are randomly chosen following the uniform distribution
$\mathcal{U}([0.8,1.2])$, and the centers of the Gaussian wells $r_i$, are chosen randomly under the
constraint that $|r_i - r_{i'}| > 4 \sigma$, \add{unless explicitely specified}.
The Kohn-Sham map is discretized using a pseudo-spectral method \cite{Trefethen2000}, and solved by
a standard eigensolver.

\subsubsection{One-dimensional case}

We set $\sigma = 0.05$ and we generated $4$ data sets using different number of wells, $n_g$, which in this case is also equal
to the number of electrons $n_e$, ranging from $2$ to $8$.

The number of discretization points is $N=320$.
We trained the architecture defined in \cref{sec:hnn} for each $n_g$, setting the number of
levels $L = 6$, using different values for $r$ and $K$.

\cref{table:results_KS_1D_no_overfitting} shows that there is no overfitting, even at this
level of accuracy and number of parameters. This behavior is found in all the numerical examples,
thus we only report the test error in what follows.

\begin{table}{}
    \begin{center}
        \begin{tabular}{cccc}
            \hline
            $r $ & \Nparams &  Training error & Validation error  \\
            \hline \hline
            2     & 2117   & 6.7e-4 & 6.7e-4 \\ \hline
            4     & 5183   & 3.3e-4 & 3.4e-4 \\ \hline
            6     & 9833   & 2.8e-4 & 2.8e-4 \\ \hline
            8     & 16067  & 3.3e-4 & 3.3e-4 \\ \hline
            10    & 33013  & 1.8e-4 & 1.9e-4 \\
            \hline \hline
        \end{tabular}
    \end{center}
    \caption{ Relative error on the approximation of the Kohn-Sham map for different $r$, with $K =
    6$, \Ntrainsample =$16000$, and \Ntestsample =$4000$.}
    \label{table:results_KS_1D_no_overfitting}
    \vspace{-.3cm}
\end{table}

\begin{figure}[ht]
    \centering
    \includegraphics[trim={10mm 5mm 5mm 8mm}, clip, width=0.95\textwidth]{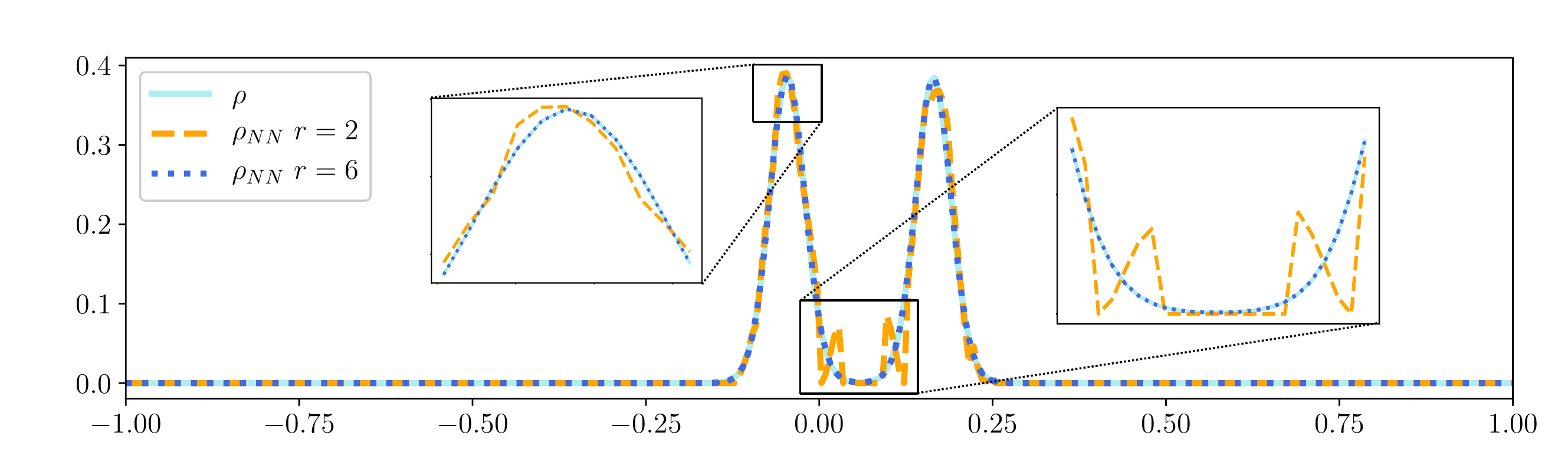}
    \caption{\label{fig:error_alpha2} Estimation using two different multiscale networks with $r = 2$, and $r =6$; with $K=6$,
    and $L = 5$ fixed.}
\end{figure}

From \cref{table:results_KS_1D_r_average} we can observe that as we increase $r$ the error decreases
sharply. \cref{fig:error_alpha2} depict this behavior. In
\cref{fig:error_alpha2} we have that if $r=2$, then the network output
$\rho_{NN}$, fails to approximate $\rho$ accurately; however, by modestly increasing $r$, the network is
able to accurately approximate $\rho$.

However, the accuracy of the network stagnates rapidly. In fact, increasing $r$ beyond $10$
does not provide any considerable gains. In addition, \cref{table:results_KS_1D_r_average}
shows that the accuracy of the network is agnostic to the number of Gaussian wells present in the
system.

\begin{table}[ht]
    \begin{center}
        \begin{tabular}{cccccccccccc}
            \hline
             $n_g\backslash r$ &  $2$ &  $4$ & $6$ &  $8$ & $10$ \\
            \hline \hline
            2    & 6.7e-4  & 3.3e-4   & 2.8e-4      & 3.3e-4   & 1.8e-4   \\  \hline
                        4    & 8.5e-4  & 4.1e-4   & 2.9e-4      & 3.8e-4   & 2.4e-4   \\  \hline
                        6    & 6.3e-4  & 4.8e-4   & 3.8e-4      & 4.0e-4   & 4.2e-4   \\  \hline
                        8    & 1.2e-3  & 5.1e-4   & 3.7e-4      & 4.5e-4   & 3.7e-4   \\
            \hline  \hline
        \end{tabular}
    \end{center}
    \caption{ Relative test error on the approximation of the Kohn-Sham map for different ranks $r$,
     with fixed $K = 6$ and \Ntrainsample =$16000$.} \label{table:results_KS_1D_r_average}
    \vspace{-.3cm}
\end{table}

In addition, we studied the relation between the quality of the approximation and $K$. We fixed
$r = 6$, and we trained several networks using different values of $K$, ranging from $2$, {\it i.e.},
a very shallow network, to $10$. The results are summarized in
\cref{table:results_KS_1D_cnn_average}. We can observe that the error decreases sharply as the depth
of the network increases, but it rapidly stagnates as $K$ becomes large.

\begin{table}{}
    \begin{center}
        \begin{tabular}{cccccc}
            \hline
            $n_g\backslash K$ &  $2$ &  $4$ & $6$ &  $8$ & $10$  \\
            \hline \hline
            2    &  1.4e-3  & 3.1e-4  &  2.8e-4  &  3.5e-4  &  2.3e-4  \\ \hline
                        4    &  1.9e-3  & 5.8e-4  &  2.8e-4  &  6.0e-4  &  7.1e-4  \\ \hline
                        6    &  2.1e-3  & 7.3e-4  &  3.8e-4  &  6.7e-4  &  6.7e-4  \\ \hline
                        8    &  2.0e-3  & 8.8e-4  &  3.7e-4  &  6.7e-4  &  6.8e-4  \\
            \hline \hline
        \end{tabular}
    \end{center}
    \caption{  Relative test error on the approximation of the Kohn-Sham map for different $K$ and
    fixed rank $r = 6$, and \Ntrainsample = $16000$.}
    \label{table:results_KS_1D_cnn_average}
    \vspace{-.3cm}
\end{table}

\add{The Kohn-Sham map is a very non-linear mapping, and we demonstrate below that the non-linear activation functions in the network play a crucial role. Consider a linear network, in which we have completely
eliminated the non-linear activation functions. This linear network can be easily implemented
by setting $K =1$, and $\phi(x) = x$. We then train the resulting linear network using the same
data as before and following the same optimization procedure. Table \ref{tab:ks_linear}
shows that the approximation error can be very large, even for relatively small $n_g$.}

\begin{table}[htb]
    \centering
    \begin{tabular}{cccc}
      \hline
      $n_g$ & \Nparams &  Training error   &  Validation error\\ \hline\hline
      $2$   &  8827    &   9.5e-2          &   9.5e-2  \\ \hline
      $4$   &  8827    &   9.7e-2          &   9.5e-2  \\ \hline
      $6$   &  8827    &   9.7e-2          &   9.7e-2  \\ \hline
      $8$   &  8827    &   8.7e-2          &   8.6e-2  \\
      \hline
      \hline
    \end{tabular}
    \caption{\label{tab:ks_linear} \add{Relative error in approximating the Kohn-Sham map for different
      number of Gaussians $n_g$ for 1D using a linear MNN network with $K=1$, $r=10$ and \Ntrainsample $=16000$.} }
\end{table}

\add{In addition, we note that the favorable behavior of MNN with respect to $n_g$ shown in the prequel is rooted in
the fact that the band gap (here the band gap is equal to $\epsilon_{n_g+1}-\epsilon_{n_g}$) remains approximately the same as we increase $n_g$. According to the density functional perturbation theory (DFPT), the Kohn-Sham map is relatively insensitive to the change of the external potential. This setup mimics an insulating system. On the other hand, we may choose $\sigma$ in \eqref{eqn:gaussian_wells} to be $0.15$, and relax the constraint between Gaussian centers to $|r_i - r_{i'}| > 2 \sigma$. The rest of the coefficients are randomly chosen using
the same distributions as before.  In this case, we generate a new data set,
in which the average gap for the generated data sets are $0.2$, $0.08$, $0.05$ and $0.03$
for $n_g$ equal to $2$, $4$, $6$, and $8$, respectively. The decrease of the band gap with respect to the increase of $n_g$ resembles the behavior of a metallic system, in which the the Kohn-Sham map becomes more sensitive to small perturbations of the potential.
After generating the samples, we trained two
different networks: the MNN with $K=6$, $r=10$ and a regular CNN with $15$ layers, $10$ channels, window size $13$. The results
are shown in Tables \ref{tab:ks_small_gap} and \ref{tab:ks_small_gap_CNN}.
From Tables \ref{tab:ks_small_gap} and \ref{tab:ks_small_gap_CNN} we can observe that MNN outperforms
CNN and that as the band gap decreases, the performance gap between MNN and CNN widens.
We point out that it is possible to partially alleviate this adverse dependence
on the band gap in the MNN by introducing a nested hierarchical structure to the interpolation
and restriction operators as shown in \cite{fan2018mnnh2}.}
\begin{table}[htb]
    \centering
    \begin{tabular}{cccc}
      \hline
      $n_g$ & \Nparams &  Training error  &  Validation error \\ \hline \hline
      2     & 20890          &  1.5e-3    &   1.9e-03          \\ \hline
      4     & 20890          &  7.1e-3    &   8.8e-03          \\ \hline
      6     & 20890          &  1.6e-2    &   2.3e-02          \\ \hline
      8     & 20890          &  2.9e-2    &   2.9e-02          \\
      \hline
      \hline
    \end{tabular}
    \caption{\label{tab:ks_small_gap} \add{Relative error in approximating the Kohn-Sham map for different
      number of Gaussians $n_g$ for 1D using an MNN network with $K=6$, $r=10$ and \Ntrainsample $=16000$.} }
\end{table}
\begin{table}[htb]
    \centering
    \begin{tabular}{cccc}
      \hline
      $n_g$ & \Nparams  &  Training error   &  Validation error  \\ \hline\hline
      2     &    19921     &  1.6e-3        &  1.9e-3  \\
      4     &    19921     &  1.7e-2        &  1.8e-2  \\
      6     &    19921     &  6.2e-2        &  6.5e-2  \\
      8     &    19921     &  9.0e-2        &  9.3e-2  \\
\hline
      \hline
    \end{tabular}
    \caption{\label{tab:ks_small_gap_CNN} \add{ Relative error in approximating the Kohn-Sham map for different
      number of Gaussians, $n_g$, for 1D using a periodic Convolutional Neural Network with $15$ layers, $10$ channels, windows size $13$ and \Ntrainsample $=16000$. }}
\end{table}

\subsubsection{Two-dimensional case}

\begin{figure}[ht]
  \centering
    \subfloat[$\rho_{NN}$]{
    \label{fig:sol_KS_2D_rho}
    \includegraphics[trim={0mm 0mm 0mm 0mm}, clip, width=0.45\textwidth]{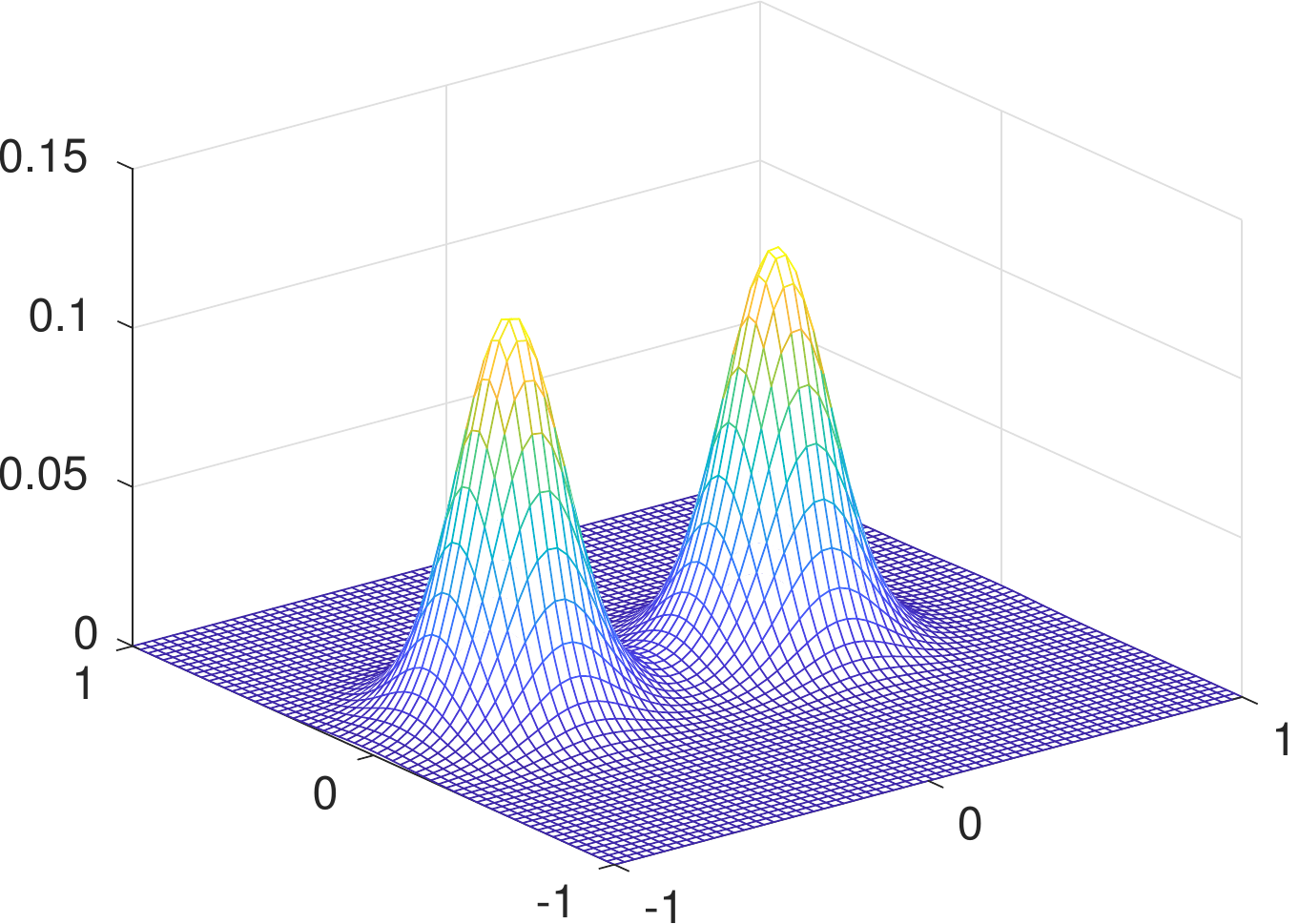}
    }
    \subfloat[$|\rho_{NN}-\rho|$]{
    \label{fig:sol_KS_2D_rho_error}
    \includegraphics[trim={0mm 0mm 0mm 0mm}, clip, width=0.45\textwidth]{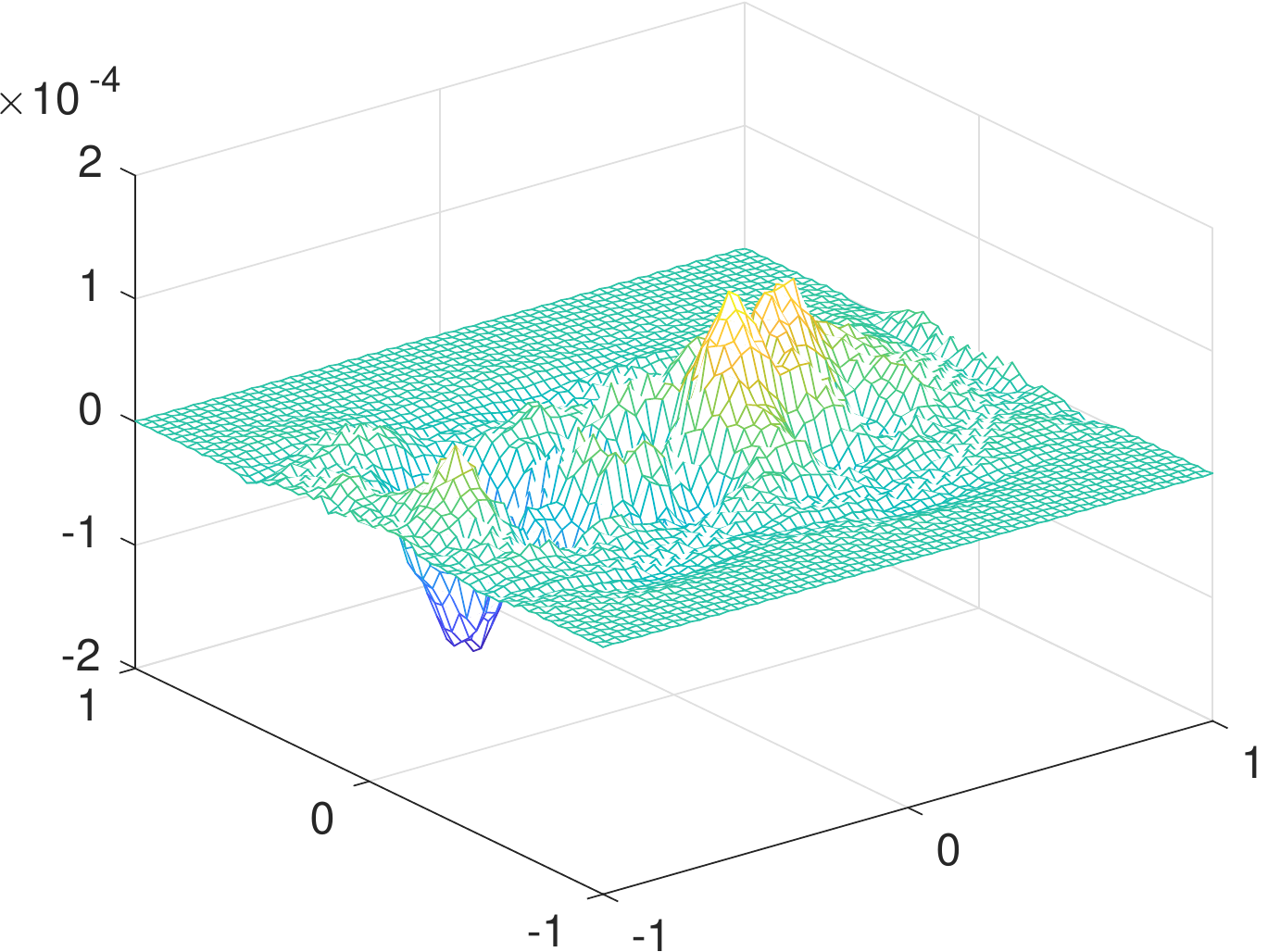}
    }
    \caption{\label{fig:sol_KS_2D} (a)  Output of the trained network on a test sample for $K=6$, and
    $\alpha = 10$; (b) error with respect to the reference solution.}
\end{figure}

The discretization is the standard extension to 2D using tensor products, using a $64\times 64$ grid.
In this case we only used $n_g=2$ and we followed the same number of training and test
samples as that in the $1$D case. We fixed $K = 6$, $L=4$, and we trained the network for different number of
channels, $r$. The results are displayed in \cref{table:results_KS_2D}, which shows the
same behavior as for the 1D case, in which the error decays sharply and then stagnates, and there is no
over fitting. In particular, the network is able to effectively approximate the Kohn-Sham map as shown in
\cref{fig:sol_KS_2D}. \cref{fig:sol_KS_2D_rho} shows the output of neural network for a test sample
and \cref{fig:sol_KS_2D_rho_error} shows the approximation error with respect to the reference.

\begin{table}{}
    \begin{center}
        \begin{tabular}{ccc}
             $r$ &  Training error & Validation error \\
             \hline \hline
             4    & 5.2e-3  & 5.2e-3    \\  \hline
             6    & 1.6e-3  & 1.7e-3    \\  \hline
             8    & 1.2e-3  & 1.1e-3      \\  \hline
             10   & 9.1e-4  & 9.3e-4    \\
             \hline \hline
        \end{tabular}
    \end{center}
    \caption{ Relative errors on the approximation of the Kohn-Sham map for 2D case for different
    $r$ and $K = 6$, \Ntrainsample =$16000$ and \Ntestsample =$4000$.} \label{table:results_KS_2D}
    \vspace{-.3cm}
\end{table}

 \section{Conclusion}\label{sec:conclusion}

We have developed a multiscale neural network (MNN) architecture for approximating nonlinear
mappings, such as those arising from the solution of integral equations (IEs) or partial differential
equations (PDEs). In order to control the number of parameters, we first rewrite the widely used
hierarchical matrix into the form of a neural network, which mainly consists of three sub-networks:
restriction network, kernel network, and interpolation network. The three sub-networks are all
linear, and correspond to the components of a singular value decomposition. We demonstrate that such
structure can be directly generalized to nonlinear problems, simply by replacing the linear kernel network
by a multilayer kernel network with nonlinear activation functions. Such ``nonlinear singular value
decomposition operation'' is performed at different spatial scales, which can be efficiently
implemented by a number of locally connected (LC) networks, or convolutional neural networks
(CNN) when the mapping is equivariant to translation. Using the parameterized nonlinear
Schr\"odinger equation and the Kohn-Sham map as examples, we find that MNN can yield accurate
approximation to such nonlinear mappings. When the mapping has $N$ degrees of freedom, the
complexity of MNN is only $O(N\log N)$. Thus the resulting MNN can be further used to accelerate the
evaluation of the mapping, especially when a large number of evaluations are needed within a certain
range of parameters.

In this work, we only provide one natural architecture of multiscale neural network based on
hierarchical matrices. The architecture can be altered depending 
on the target application. Some of the possible modifications and extensions are listed below. 
1) In this work, the neural network is inspired by a hierarchical matrix with a special case of strong 
admissible condition. One can directly construct architectures for $\cH$-matrices with the weak
admissible condition, as well as other structures such as the fast multiple methods, $\cH^{2}$-matrices and wavelets. 
2) The input, $u$, and output, $v$, in this work are periodic. The network can be directly extended 
to the non-periodic case, by replacing the periodic padding in $\LCK$ by some other padding
functions. One may also explore the mixed usage of LC networks and CNNs in different components of
the architecture.
3) The matrices $A\sps{\ad}$ and $M\sps{\ell}$ can be block-partitioned in different
ways, which would result in different setups of parameters in the $\LCK$ layers.
4) The $\LCR$ and $\LCI$ networks in Algorithm \ref{alg:mnn} can involve nonlinear activation functions 
as well and can be extended to networks with more than one layer.
5) The $\LCK$ network in Algorithm \ref{alg:mnn} can be replaced by other architectures. 
In principle, for each scale, these $\LCK$ layers can be altered to any network, for example 
the sum of two parallel subnetworks, or the ResNet structure\cite{he2016deep}. 6) It is known that
$\cH$-matrices can well approximate smooth kernels but become less efficient for highly oscillatory
kernels, such as those arising from the Helmholtz operator in the high frequency regime. The
range of applicability of the MNN remains to be studied both theoretically and numerically.

\section*{Acknowledgements}
The authors thank Yuehaw Khoo for constructive discussions.
The work of Y.F. and L.Y. is partially supported by the U.S. Department of Energy, Office of
Science, Office of Advanced Scientific Computing Research, Scientific Discovery through Advanced
Computing (SciDAC) program and the National Science Foundation under award DMS-1818449,
and the GCP Research Credits Program from Google. The work of Y.F. is also partially supported by
AWS Cloud Credits for Research program from Amazon.
The work of L.L and L. Z. is partially supported by the Department of Energy under Grant No.
DE-SC0017867,  by the Department of Energy under the CAMERA project, and by the Air Force Office of Scientific Research under award number FA9550-18-1-0095.

\bibliographystyle{abbrv}
\bibliography{../nn}
\end{document}